\documentclass{amsart}
\usepackage{amsmath}
\usepackage{amsthm}
\usepackage{amssymb}
\usepackage{hyperref}

\newtheorem{thm}{Theorem}

\newtheorem{prop}{Proposition}
\newtheorem{rmk}{Remark}
\newcommand{\be}{\begin{equation}}
 \newcommand{\ee}{\end{equation}}
 \newcommand{\ba}{\begin{array}{l}}
 \newcommand{\ea}{\end{array}}
 \newcommand{\pa}{\partial}
\newcommand{\na}{\nabla}

 \newcommand{\la}{\label}
\newcommand{\fr}{\frac}

\newcommand{\dm}{{\mbox{div}}_m}

\newcommand{\nam}{\nabla_m}
\newcommand{\dam}{\mbox{div}_m}

\newcommand{\nax}{\nabla_x}
\newcommand{\naa}{\nabla_a}
\newcommand{\naxi}{\nabla_{\xi}}
\newcommand{\wf}{\widehat{f}}
\newcommand{\og}{\overline{\gamma}}
\newcommand{\Dx}{\Delta_x}
\newcommand{\dx}{{\mbox{div}}_x}
\newcommand{\ws}{\widetilde{\sigma}}

\newcommand{\Rr}{{\mathbb{R}}}
\newcommand{\e}{\epsilon}

\begin{document}

\title
{Remarks on Oldroyd-B and Related Complex Fluid Models}

\author{Peter Constantin\and Weiran Sun}

\address{Department of Mathematics, The University of Chicago\\5743 S. University Ave., Chicago, Illinois 60637}
\email{const@math.uchicago.edu, wrsun@math.uchicago.edu}



\begin{abstract}
We prove global existence and uniqueness of solutions of Oldroyd-B systems
with relatively small data in  $\Rr^d$, in a large functional setting ($C^{\alpha}\cap L^1$). This is a stability result, solutions select an equilibrium and converge exponentially to it. Large spatial derivatives of the initial density and stress are allowed, provided the $L^{\infty}$ norm of the density and stress are small enough. We prove global regularity for large data for a model in which the potential responds to high rates of strain in the fluid. We also prove global existence for a class of large data for a didactic scalar model which attempts to capture in the simplest way the essence of the dissipative nature of the coupling to fluid. This latter model has an unexpected cone invariance in function space that is crucial for the global existence.
\end{abstract}

\maketitle

\section{Introduction}
The complex fluids models we consider treat the interaction of a normal fluid (incompressible in our case) with insoluble matter suspended in it. Models have been devised to deal with microscopic elastic thread-like objects such as polymers (\cite{doied}, \cite{otbook}). The complicated hydrodynamic interactions are simplified using the separation of scales, replacing the many degrees of freedom due to them by few representative ones. In the models we consider here, an end-to-end vector in $\Rr^3$ ($m$ in this paper) represents the orientation of the polymer and it is the sole variable retained to describe the microscopic object. The models that are most studied with this degree of simplification include the kinetic description of the Oldroyd B model, the FENE model and variants. The probability distribution of the vectors $m$ depends on time and physical space, and it is assumed to be absolutely continuous with respect to the usual Lebesgue measure $dm$, so the chance of finding $m$  at time $t$ and location $x\in \Rr^3$ in the volume element $dm$ is $f(x,m,t)dm$. The aim is to describe the evolution of $f$ in space and time. When the polymer concentration is so dilute that
the polymers do not interact, and when the fluid is smooth, then the equation of $f$ is a linear Fokker-Planck equation. The effect of the fluid can be understood perturbatively: because the objects suspended are assumed to have faster time scales and shorter length scales than the scales on which solvent varies, it is then 
justified to treat the fluid as varying little as far as the microscopic suspended objects are concerned. This {\em{macro-micro}} interaction can be rationally discussed, and leads to a  kinetic description of the particles, given a fluid flow. On the other hand, the collective effect that the particles have on the flow itself is a macroscopic effect of microscopic insertions: a {\em{micro-macro}} interaction. These interactions are much more mysterious. There have been attempts to produce systematic upscalings based on non-equilibrium thermodynamics formalisms (\cite{beris}, \cite{grot}). When the microscopic insertions have a larger number of degrees of freedom, $m\in M$, where $M$ is a Riemannian manifold representing a finite number of degrees of freedom with constrains, then the kinetic equation is a Smoluchowski equation on the manifold. The upscaling principle advocated by one of us (\cite{c}, \cite{c-smo}) is easiest formulated as the {\em{requirement}} that the sum of the energy of the fluid and the free energy of the microscopic quantities be a Lyapunov functional for
the coupled system. In the known examples, this requirement leads to familiar rules of determining the added polymeric stress from the micro-micro and the macro-micro interactions (\cite{lebris}).
The mathematical treatment of the coupled systems is far from being complete but has received recently more attention. Early work (\cite{saut}, \cite{renardy}) established local existence results for Oldroyd-B and FENE type equations. The problem of existence of weak solutions is rather open. Global existence of weak solutions in the presence of spatial diffusion of the polymers was proved in a sequence of papers, the most recent of which is (\cite{suli}). Global existence of weak solutions via propagation of compactness was proved under the corotational assumption (\cite{lionsmasmoudi}, \cite{lionsmasmoudi1}) and, very recently, for the full FENE model (\cite{mas1}). There is no such result for the Oldroyd B model. The global existence of smooth solutions for small data for Oldroyd B-type models was established in (\cite{leiz}, \cite{fanghuasmall}).
Global existence of smooth solutions for large data in 2D was established for Smoluchowski equations on compact manifolds (\cite{c-smo}, \cite{cftz}, \cite{cm}, \cite{cs}, \cite{cs1}, \cite{otto}). Global regularity for large data in the FENE case, under the corotational assumption was proved in (\cite{fanghua}, \cite{mas}). An approach based on Lagrangian particle dynamics was described in (\cite{fanghuachun}). Sufficient conditions for regularity in terms of bounds on the added stress tensor were established in (\cite{chemm}, \cite{kmt}) and further refined in (\cite{leim}). Numerical evidence for singularities was provided in (\cite{thomasses}).

In this paper we address issues related to global existence of smooth solutions in simplest kinetic models. In the presence of a quiescent solvent, the polymer distribution is the unique stable time independent solution of the linear Fokker-Planck equation
$$
\pa_t f = \epsilon\dm\left(f\nam (\log f + U)\right )  
$$
i.e. $f=Z^{-1}e^{-U}$ with $Z$ a normalizing constant. The simplest form of $U$ is harmonic, $U(m) = \fr{|m|^2}{2R^2}$ and then the solution is Gaussian. The constant $R^2$ represents the expected value $<|m|^2>$ of the square of the length of the end-to-end vector. In the presence of smooth fluids the kinetic equation changes and acquires a dependence on the macroscopic independent variable $x$. 

In the second section of the paper we provide a priori bounds for linear Fokker-Planck equations with Oldroyd B type potentials. We do this by essentially giving a formula for the solution of the Fokker-Planck equation in terms of the  particle trajectories of the underlying flow.

In the third section we derive a classical estimate for the velocity gradient in terms of the added polymeric stress, when the coupling to fluid is done via time-independent Stokes equations. 

In the fourth section we prove global existence results for small data in $C^{\alpha}$, at arbitrary Deborah numbers. The previously known global existence results for small data are done using energy estimates and require derivatives of the stress. The main difficulty in obtaining bounds for small data is the fact that quantities are not  monotonically decreasing, and, in particular, the spatial gradient of the density of particles can grow in time, but it nevertheless saturates. The system formed by the kinetic equation for the particle distribution and the velocity of the fluid can be reduced to a closed system for a reduced added stress $\tau$ and the particle density. This system looks like a damped and driven Euler equation where $\tau$ (a matrix) plays the role of vorticity and $\rho$, the transported density is part of the driver. The global existence of solutions is proved under the assumption that the $L^{\infty}$ norm of the initial data is small enough
(\ref{smeps}). The class of initial data that lead to global existence includes densities and stresses with arbitrarily large gradients, provided the $L^{\infty}$ norms are suitably small. While the limit added stress $\tau$ and limit velocity vanish, the limit particle density is a re-arranged initial particle density, a deformation of the initial density along the limit back-to-labels map. The class of solutions we discuss is in fact rather wide, and in particular the gradients need to be understood in distribution sense. The proof of uniqueness of solutions is done using Lagrangian transformations. An Eulerian proof is not readily available at this moment. In order to explain this we pursue the analogy with Euler equations. Both the damped Euler vorticity equation and the equation for  $\tau$ have the form
$$
D_t \tau = -c\tau + N(\tau)
$$
where $c\ge 0$ is the damping and the nonlinearity $N$ is quadratic and non-local, obtained from products of $\tau$ and principal value singular integrals ${\mathbb H}\tau$.
Here $D_t$ is material derivative along the divergence-free velocity $u$.
The local (and global for small data) existence results are in Banach algebras in which $N$ is a continuously differentiable nonlinear function. One of the largest and most natural such Banach algebras is $C^{\alpha}$, $0<\alpha<1$. 
If we consider the difference $\delta$ of two solutions $\tau_1$, $\tau_2$, it obeys the equation
$$
{\overline D}_t \delta = - c\delta + N'(\overline{\tau})\delta - v\cdot\nax{\overline{\tau}}
$$
where ${\overline D}_t$ is transport along the average velocity ${\overline u}= \fr{1}{2}(u_1 + u_2)$, $\overline{\tau}$ is the average of the solutions,
${\overline{\tau}} = \fr{1}{2}(\tau_1 + \tau_2)$, $N'(\overline{\tau})$ is the linearization of $N$ obtained using polarization, and $v$ is the velocity difference $v= u_1-u_2$. The last term in the equation is present because the velocities depend on $\tau$: $u_1$ is obtained from $\tau_1$ and $u_2$ from $\tau_2$. In both the (damped-driven) Euler case and in the Stokes-Fokker-Planck system we study, the dependence of $u$ on $\tau$ is linear. The term $v$ depends therefore linearly on $\delta$.  Uniqueness proofs based on energy estimates proceed by estimating the growth of $\delta$ in some norm. As the initial data vanishes, any closed, super-linear estimate of the  norm is sufficient to prove uniqueness. (We could allow logarithmically sublinear estimates, but that is as far as that can go, in general). There are two sources of difficulty in obtaining these estimates. The transport term ${\overline u}\cdot\nax\delta$ can be dealt by using integration by parts, if the norm of $\delta$ we consider is an $L^p$ norm.
The other source of difficulty is $\nax{\overline{\tau}}$. If we are in the framework of $C^{\alpha}\cap L^1$ adopted in this paper, then $v$ is very nice but
in order to deal with $\nax\overline \tau$ we need to work in $H^{-1}$. In general, the presence of both $v\cdot\nax\overline\tau$  and $\overline{u}\cdot\nax\delta$ together is lethal for an Eulerian, energy method approach, unless we work in a space in which $\nax\overline\tau$ is tame, like for instance $W^{1,p}(\Rr^d)$, with $p>d$.
The reason there is a successful  Eulerian proof of uniqueness for Euler equations for the vorticity in $C^{\alpha}$ is that, for Euler equations, the vorticity equation comes from a velocity equation with a good cancellation property. Thus we can ``retreat'' to one less derivative in the 
equation for $\delta$. In terms of the vorticity equation this means integrating against a vectorial stream function, i.e. working in $H^{-1}$. In the Euler equtions the transport term together with the term $N'({\overline{\tau}})\delta$ can be integrated by parts one more time, revealing only $H^{-1}$ terms in $\delta$. In the absence of such algebraic reasons, the retreat to one less derivative does not appear to work, and the fix for one difficult term does not work for the other, and vice-versa. The Lagrangian approach introduces commutators that are well behaved, and is successful.

In the fifth section we show that the system is regularized for large data if the potential is allowed to respond to excessively high rates of strain $S$ in the fluid. The Oldroyd-B potential is harmonic, $U= \fr{|m|^2}{2R^2}$, as it were generated by a spring force that is related to a material property of the polymers. But, in fact, the representation of the polymer by $m$ is over-simplified, and $R$ has something to do with an average restoring force in the ensemble, and not with the maximum extension allowed for an individual molecule. In the case there is feedback coupling to a fluid, the interaction might have properties that depend of the local properties of the fluid. We show that, if we allow $R$ to depend on the local rate of strain in the fluid, when this exceeds a threshold, and to grow with it, $\fr{D_t R}{R} = \delta(S)$, then the equation has global smooth solutions for all large data. Here $\delta(S)$ can vanish if the rate of strain $S$ is not too high, and becomes asymptotically linear in $S$ for high values of (norms of) $S$. Allowing $R$ to grow in response to high rate of strain is like allowing more entropic ``slack'' in the molecules at high rates of strain, and this turns out to be a mathematical regularizing mechanism.

In the sixth section we return to the Oldroyd B large data problem, and give an example of blow up in one dimension for a crude model that does not respect the principle of free energy  decay. We motivate and give then a simple scalar example of an equation that has some of the features of the type of equations arising in the large data Oldroyd B equation, and which has an unbounded set of initial data that lead to global smooth behavior. The reason for the global existence is an invariance of a cone in function space under the nonlinear evolution. The consequnece of this invariance is an a priori bound that is sufficient for the persistence of smoothness of solutions.

\section{Linear Fokker-Planck Equation: kinematic observations}
We take a vector field $u(x,t)$ and a scalar $f(x,m,t)$ representing a two or three-dimensional incompressible velocity and a particle distribution.
We start by describing the simplest particle distribution equation
\be
D_tf +((\nax u)m)\cdot \nam f = \e\dam(f\nam(\log f + U(m)))
\la{lfp}
\ee
Here $D_t = \pa_t + u\cdot\nax$ is material derivative.
The potential is given by
\be
U(m) = \frac{|m|^2}{2R^2}.
\la{um}
\ee 
In this section $R$ is a positive constant. We consider here the case $m\in\Rr^d$.
We associate to the particle distribution $f$ an added stress tensor:
\be
\sigma (x,t) = \int_{m\in\Rr^d} \left(m\otimes \nam U\right )f(x,m,t)dm
\la{sigmam}
\ee
In this normalization $R$ has units of length, $\e/R^2$ is an inverse time. The potential $U$ and the stress $\sigma$ are nondimensional.  We study first (\ref{lfp}), without regard to the coupling to the fluid. In what follows, the velocity field could be quite arbitrary, but it is assumed to be smooth enough for the calculations below: the variable $x\in \Rr^d$, and the velocity field $u(x,t)$ is in $C^{1,\alpha}(\Rr^d)$, is divergence-free
$$
\dx u=0
$$
and decays at infinity. 
Before we start our computations, let us make a few general comments regarding the matrix $\sigma$. First, because the potential $U$ is radially symmetric (\ref{um}) and from its definition (\ref{sigmam}), it follows that
$\sigma $ is symmetric and non-negative. Moreover its off-diagonal entries are bounded by the trace
\be
|\sigma^{ij}(x,t)| \le\frac{1}{2}Tr(\sigma(x,t)).
\la{indi}
\ee
We take now the Fourier transform in $m$
\be
\widehat{f}(x,t, \xi) = \int_{\Rr^d}e^{-im\cdot\xi}f(x,m,t)dm
\la{fouf}
\ee
In view of (\ref{um}), $\widehat{f}$ obeys
\be
D_t\widehat{f}(x,\xi,t) + \left[\frac{\epsilon}{R^2}{\mathbb I} - (\nax u)^T\right ]\xi\cdot\naxi\wf = -\e|\xi|^2\wf(x,\xi,t)
\la{lfpfou}
\ee
We used the fact that $u$ is divergence-free and denoted $(\nax u)^T$ the transposed matrix. We solve (\ref{lfpfou}) on characteristics. 
The connection between the Lagrangian paths of $u$ and the Oldroyd B equation was exploited in \cite{fanghuachun} in the less complicated case $\e=0$, $\rho =1$. Let $X(a,t)$ be the particle paths,
\be
\pa_tX(a,t) = u(X(a,t),t)
\la{X}
\ee
with $X(a,0) =a$ and let
\be
g(a,t)= (\nax u)(X(a,t),t)\la{g}
\ee
Let now $\xi = \xi(a,\eta,t)$ solve the ODE
\be
\frac{d}{dt} \xi = \frac{\e}{R^2}\xi - g(a,t)^T\xi
\la{xieq}
\ee
with initial data $\xi(a, \eta, 0) = \eta$.  We take the fundamental matrix 
$\Phi(a,t)$, solution of the linear ODE system
\be
\frac{d}{dt}\Phi(a,t) = - g(a,t)^T\Phi(a,t)
\la{phiq}
\ee
with initial data $\Phi(a,0) = \mathbb I$ and we have then
\be 
\xi(a,\eta, t) = e^{\frac{\epsilon t}{R^2}}\Phi(a,t)\eta
\la{xiph}
\ee
We write now
\be 
F(a,\eta,t) = \wf(X(a,t),\xi(a,\eta,t), t)\la{F}.
\ee
Then (\ref{lfpfou}) implies
$$
\frac{d}{dt} F(a,\eta, t) = - \e |\xi(a,\eta, t)|^2 F(a, \eta, t)
$$
and integrating we obtain
\be
F(a,\eta, t) = e^{-\epsilon \int_0^t|\xi(a,\eta, s)|^2ds}\wf(a,\eta,0)
\la{Feq}
\ee
Now we invert the linear map $\eta \mapsto \xi(a,\eta,t)$, and write from (\ref{xiph}), 
\be
\eta(a,\xi, t) = e^{-\fr{\e t}{R^2}}\Psi(a,t)\xi
\la{etat}
\ee
where $\Psi(a,t) = \Phi(a,t)^{-1}$ obeys
\be
\frac{d}{dt}\Psi(a,t) = \Psi(a,t) g(a,t)^T.
\la{psieq}
\ee
with initial data $\Psi(a,0) = \mathbb I$. Reading (\ref{F}) at $\eta = \eta(a,\xi,t)$ for a fixed $\xi$ and substituting in (\ref{Feq}) we obtain
\be
\wf(X(a,t), \xi, t) = e^{-\e\int_0^t|\xi (a, \eta(a, \xi, t), s)|^2ds}\wf_0(a, e^{-\frac{\e t}{R^2}}\Psi(a,t)\xi)\la{wfa}
\ee
where $\wf_0(x,\xi)$ is the Fourier transform in $m$ of the initial data $f_0(x,m) = f(x,m,0)$.
Let us consider now the ``back-to-labels'' map $A(x,t)$, inverse of $X(a,t)$.
Let 
\be 
M(x,t) = \Psi(A(x,t),t)\la{M}
\ee
and note that it obeys the transport equation
\be 
D_tM = M (\nax u)^T
\la{Meq}
\ee
with initial data $M(x,0) = \mathbb I$. Reading (\ref{wfa}) at $a= A(x,t)$ and we deduce
\be
\wf( x,\xi, t) = e^{-\e \int_0^t|\xi(A(x,t), \eta(A(x,t),\xi,t), s)|^2ds}\wf_0(A(x,t), e^{-\frac{\e t}{R^2}}M(x,t)\xi).
\la{wfxprlim}
\ee
We compute directly from (\ref{xiph}) and (\ref{etat})
\be
\xi(A(x,t), \eta(A(x,t), \xi, t), s) = e^{-\fr{\e (t-s)}{R^2}}Q(x,t,s)\xi
\la{xieta}
\ee
with
\be
Q(x,t,s) = q(A(x,t), t, s), \quad {\mbox{for}}\; \;t\ge s,\la{Qq}
\ee
\be
q(a,t,s) = \Phi(a,s)\Psi(a,t)   \quad {\mbox{for}}\; \;t\ge s, \la{q}
\ee
and deduce thus from (\ref{wfxprlim})
\be
\wf(x,\xi,t) = e^{-\e\int_0^te^{-\fr{2\e (t-s)}{R^2}}\left|Q(x,t,s)\xi\right|^2ds}\wf_0(A(x,t), e^{-\fr{\e t}{R^2}}M(x,t)\xi)
\la{wfx}
\ee
Let us note now that the matrix $\sigma^{ij}(x,t)$ is computed from the Hessian of $\wf$ at $\xi=0$,
\be
\sigma^{ij}(x,t) = -\frac{1}{R^2}\fr{\pa^2\wf }{\pa \xi_i\pa \xi_j}(x,\xi,t)_{\left | \xi =0\right.}
\la{sigfou}
\ee
Using (\ref{wfx}) and noting that the cross-terms vanish, we deduce
\be
\ba
\sigma^{ij}(x,t) = \frac{2\e }{R^2}\rho(x,t)\int_0^te^{-\fr{2\e (t-s)}{R^2}}\left[Q^T(x,t,s)Q(x,t,s)\right]_{ij}ds \\ - \frac{e^{-\frac{2\e t}{R^2}}}{R^2}M_{ki}(x,t)M_{lj}(x,t)\fr{\pa^2\wf_0 }{\pa \xi_k\pa \xi_l}(A(x,t),\xi)_{\left | \xi =0\right.}
\ea
\la{sigfo}
\ee 
We denoted by $\rho(x,t) = \wf(x,0,t) = \int_{\Rr^d} f(x,m,t)dm$. In view of (\ref{lfp}) this obeys the transport equation
\be
D_t\rho = 0\la{dtrho}
\ee
and therefore it is given by
\be
\rho(x,t) =\wf_0(A(x,t))= \rho_0(A(x,t))\la{rhorho0}
\ee
in terms of the initial particle density at $x$. We used this in (\ref{sigfo}) as well as the summation convention. In terms of the initial stress, the expression is 
\be
\ba
\sigma^{ij}(x,t) =   \frac{2\e }{R^2}\rho(x,t)\int_0^te^{-\fr{2\e (t-s)}{R^2}}\left[Q^T(x,t,s)Q(x,t,s)\right]_{ij}ds\\
+ {e^{-\frac{2\e t}{R^2}}}M_{ki}(x,t)M_{lj}(x,t)\sigma^{kl}(A(x,t), 0).
\ea
\la{sigsig0}
\ee
The solution in (\ref{sigsig0}) solves the equation 
\be
D_t\sigma = -\frac{2\e}{R^2}\sigma + (\nax u)\sigma + \sigma (\nax u)^T + \fr{2\e}{R^2}\rho(x,t)\mathbb I
\la{sigeq}
\ee
which can be easily derived from the equation (\ref{lfp}) by multiplying with $\frac{1}{R^2}m_im_j$ and integrating $dm$. 
Let $\theta(x,t)$ be a passive scalar, i.e., a solution of
\be
D_t\theta = 0.
\la{thetaeq}
\ee
Then, if $M$ solves (\ref{Meq}) then $M\nax\theta$ is again a passive scalar, i.e.,
\be
D_t\left (M\nax\theta\right ) = 0.
\la{ertel}
\ee
This can be easily checked because
$$
D_t(\partial_j\theta) = - (\partial_ju_l)\partial_l\theta
$$
and so
$$
D_t(M_{ij}\partial_j\theta) = M_{ik}(\partial_k u_j)\partial_j\theta -M_{ij}(\partial_ju_l)\partial_l\theta =0
$$
This means that if $\theta(x,t) = \theta_0(A(x,t))$ with arbitrary smooth $\theta_0(a)$ then $M(x,t)\nax\theta(x,t) = (\na_a\theta_0)(A(x,t))$ holds. In fact, because the initial datum of $\Psi$ is the identity matrix, and because of uniqueness of ODEs, it follows from (\ref{psieq}) that
\be
\Psi (a,t) = (\naa X)^T(a,t)
\la{psina}
\ee
and consequently 
\be 
M_{ki}(x,t) = \frac{\partial X^i}{\partial a_k}(A(x,t),t).
\la{maq}
\ee
From (\ref{sigsig0}) and (\ref{maq}) it follows that that
\be
\ba
\sigma(x,t) = \frac{2\e }{R^2}\rho(x,t)\int_0^te^{-\fr{2\e (t-s)}{R^2}}Q^T(x,t,s)Q(x,t,s)ds\\
+ e^{-\frac{2\e t}{R^2}}(\naa X(A(x,t),t)) \sigma_0(A(x,t)) (\naa X(A(x,t),t))^T
\ea
\la{signaa}
\ee
Let us introduce the notations
\be
e^i(x,t) := (\na X^{i})(A(x,t),t)
\la{ei}
\ee
\be
e^{ij}(x,t):= \sigma_0(A(x,t))e^i(x,t)\cdot e^j(x,t)
\la{eij}
\ee
Note that
\be
\sigma (x,t) = \fr{2\e}{R^2}\rho(x,t)\int_0^te^{-\fr{2\e(t-s)}{R^2}}Q^T(x,t,s)Q(x,t,s)ds 
+ e^{-\fr{2\e t}{R^2}} \left(e^{ij}(x,t)\right)_{ij}
\la{al}
\ee
Moreover,
\be
D_t e^{i} = (\nax u)e^{i}
\la{eieq}
\ee
with initial data that are constant in space and equal the canonical basis of $\Rr^d$, $e^{i}(x,0) =(\delta_{ij})_j$ and that $q$ given in (\ref{q}) is given in terms of the gradient $\naa X$ by
\be
q(a,t,s) =[(\naa X(a,s))^{-1}]^T(\naa X(a,t))^{T}\la{qa}
\ee
solving 
\be
\pa_s q(a,t,s) = -g^T(a,s)q(a,t,s)
\la{qaeq}
\ee
and
\be
q(a, t, t) = {\mathbb I}, \quad q(a, t, 0) = (\naa X(a,t))^T.
\la{qf}
\ee
Passing to Lagrangian variables in (\ref{al}) we obtain
\be
\ba
\sigma(X(a,t),t) = \fr{2\e}{R^2}\rho_0(a)\int_0^te^{-\fr{2\e (t-s)}{R^2}}q^T(a,t,s)q(a,t,s)ds\\
+ e^{-\fr{2\e t}{R^2}}(\naa X(a,t))\sigma_0(a)(\naa X(a,t))^T
\ea
\la{all} 
\ee
and integrating by parts we deduce
\be
\ba
\sigma(X(a,t),t) -\rho_0(a){\mathbb I} = \\ 2\rho_0(a)\int_0^te^{-\fr{2\e(t-s)}{R^2}}q^T(a,t,s)S(X(a,s),s)q(a,t,s)ds + e^{-\fr{2\e t}{R^2}}(\pa_a X(a,t))\tau_0(a)(\pa_aX(a,t))^T
\ea
\la{ala}
\ee
where we introduced the reduced stress
\be
\tau (x,t) = \sigma(x,t) -\rho(x,t)\mathbb I
\la{taudef}
\ee
and where
\be
S(x,t) = \fr{1}{2}\left[(\nax u(x,t)) + (\nax u (x,t))^T\right]
\la{S}
\ee
is the rate of strain. Returning to Eulerian variables, we have
\be
\ba
\tau(x,t) =\\
2\rho(x,t)\int_0^te^{-\fr{2\e(t-s)}{R^2}}Q^T(x, t, s)S(X(A(x,t),s),s)Q(x,t,s)ds \\+ e^{-\fr{2\e t}{R^2}}(\naa X(A(x,t),t))\tau_0(A(x,t))(\naa X(A(x,t),t))^T
\ea
\la{tau}
\ee 
It is clear that in general, once the velocity $u(x,t)$ is given, we can compute everything we need to know about $\tau $ (and hence $\sigma$) from (\ref{tau}). We will be interested in particular in the relationship between various norms of $u$ and norms of $\tau$, as these relationships will serve in establishing bounds for solutions of nonlinear equations. It is clear that we should start with the fields $e_i$, and in order to understand their norms we perform some standard Lagrangian estimates.
We take the finite difference
\be
\delta_h\left(\partial_aX^i(a,t)\right) = \partial_aX^i(a+h,t)-\partial_aX^i(a,t)
\la{delx}
\ee
Its equation follows from (\ref{X}) by differentiation,
\be
\pa_t(\na_a X^i) = g(a,t)(\na_aX^i)
\la{naxeq}
\ee
and then taking the finite difference:
\be
\fr{d}{dt}\delta_h\left(\partial_aX^i(a,t)\right) = \fr{g(a,t)+g(a+h,t)}{2}\delta_h\left(\partial_aX^i(a,t)\right) + \delta_h g(a,t)\fr{\partial_aX^i(a+h,t) + \partial_aX^i(a,t)}{2}
\la{delxeq}
\ee
where $g(a,t)$ is the velocity gradient in Lagrangian coordinates (\ref{g})and
\be
\delta_h g(a,t) = g(a+h,t)-g(a,t).
\la{delg}
\ee 
For matrices $L$ we use the notation $|L|$ for the Euclidean norm of the matrix $|L| =\sqrt{Tr{L^*L}}$. We denote
\be
\gamma(t) = \sup_a |g(a,s)| = \|\nax u(\cdot, t)\|_{L^{\infty}(dx)}
\la{gamma}
\ee
a quantity of some importance in the sequel. From (\ref{naxeq}) and the fact that the initial data for $\pa_aX^i $ is $(\delta_{ij})_{j=1,\dots d}$ we have that
\be
\sup_a |\pa_aX^i(a,t)|\le e^{\int_0^t\gamma(s)ds}
\la{naab}
\ee
Because the initial data for $\delta_h\left(\partial_aX^i\right )$ vanishes, we have, from Gronwall's inequality and (\ref{naab}):
\be
\left |\delta_h\left(\pa_aX^i(a,t)\right)\right| \le e^{\int_0^t\gamma(s)ds}\int_0^t|\delta_h g(a,s)|ds
\la{grown}
\ee
We are interested in quantities in their Eulerian form. We consider the H\"{o}lder seminorm
\be
[\phi]_{\alpha} = \sup_{x\neq y}\frac{|\phi(x)-\phi(y)|}{|x-y|^{\alpha}}
\la{semin}
\ee
with $0<\alpha<1$. We note that the back-to-labels maps are Lipschitz and
\be
\sup_{x\neq y}\fr{|A(x,t)-A(y,t)|}{|x-y|}\le e^{\int_0^t\gamma(s)ds}
\la{alip}
\ee
because
\be
D_tA = 0\la{aeq}
\ee
with initial data $A(x,0) =x$, and consequently
\be
D_t\nax A = -(\nax A)(\nax u).
\la{naxaeq}
\ee
In view of (\ref{alip}) and the fact that H\"{o}lder seminorms behave nicely with respect to compositions with Lipschitz functions,
$$
[\phi\circ A]_{\alpha}\le [\phi]_{\alpha} \lambda^{\alpha}
$$
if $\lambda$ is a bound on the Lipschitz seminorm of $A$, it follows that
\be
[e^{i}(\cdot,t)]_{\alpha}\le [\na X^{i}(\cdot,t)]_{\alpha} e^{\alpha\int_0^t\gamma(s)ds}
\la{eixi}
\ee
and, similarly
\be
[g(\cdot,t)]_{\alpha}\le [\nax u(\cdot,t)]_{\alpha}e^{\alpha\int_0^t\gamma(s)ds}
\la{gnaxu}
\ee
From these considerations and from (\ref{grown}) it follows that
\be
[e^i(\cdot, t)]_{\alpha} \le e^{(1+2\alpha)\int_0^t\gamma(s)ds}\int_0^t[\nax u(\cdot, s)]_{\alpha}ds.
\la{einaxu}
\ee
The reason for the $2\alpha$ ``loss'' in comparison to the Lagrangian estimate (\ref{grown}) is that we ``loose'' twice by composing with Lipschitz functions, once because the objective of the estimate is composed, and another, because what what we estimate it with is also composed. 
We also obtain, in the same manner
\be
\|\nax e^{i}(\cdot,t)\|_{L^p(\Rr^d)} \le e^{3\int_0^t\gamma(s)ds}\int_0^t\|\nax\nax u(\cdot,s)\|_{L^p(\Rr^d)}ds
\la{naxeip}
\ee
We perform similar calculations for the matrix $Q(x,t,s)$ taking into account  (\ref{Qq}, \ref{q}, \ref{qaeq}, \ref{qf}). We deduce
\be
\|Q(\cdot, t,s)\|_{L^{\infty}(\Rr^d)} \le Ce^{\int_s^t\gamma(z)dz}
\la{qinf}
\ee
\be
|\delta_h q(a,t,s)| \le Ce^{\int_s^t\gamma(z)dz}\int_s^t|\delta_h g(a,z)|dz\la{delhq}
\ee
and consequently
\be
\left [Q(\cdot, t,s)\right]_{\alpha} \le Ce^{(1+2\alpha)\int_0^t\gamma(z)dz}
\int_s^t[\nax u(\cdot, z)]_{\alpha}dz.
\la{qalpha}
\ee
with $C =\sqrt{d}$ depending only on the dimension of space $d$. 
We are ready to state the estimates on $\tau$ in terms of $u$:
\begin{prop}\la{pasprop} Let $u$ be a divergence-free function belonging to $L^1(0,T; C^{1,\alpha}(\Rr^d)$. Let $\sigma_0(a)$ be a H\"{o}lder continuous, $L^{1}(\Rr^d)$, positive symmetric matrix and let $\rho_0$ be a positive, H\"{o}lder continous, $L^{1}(\Rr^d)$ function. Then $\tau = \sigma -\rho {\mathbb I}$, given in the 
expression (\ref{tau}), obeys
\be
\ba
\|\tau (\cdot, t)\|_{L^{\infty}(\Rr^d)} \le C\|\rho_0\|_{L^{\infty}(\Rr^d)}\int_0^t \exp{\left\{-\fr{2\e(t-s)}{R^2} + 2\int_s^t\gamma(z)dz\right\}} \gamma(s)ds\\
+C\|\tau_0\|_{L^{\infty}(\Rr^d)}\exp{\left \{-\fr{2 \e t}{R^2} + 2\int_0^t\gamma(s)ds\right\}},
\ea
\la{wsufty}
\ee
and
\be
\ba
[\tau(\cdot, t)]_{\alpha} \le C[\rho_0]_{\alpha}e^{\alpha\int_0^t\gamma(z)dz}\int_0^t \exp{\left\{-\fr{2\e(t-s)}{R^2} + 2\int_s^t\gamma(z)dz\right\}} \gamma(s)ds\\
+C\|\rho_0\|_{L^{\infty}(\Rr^d)}\int_0^t\exp{\left\{-\fr{2\e(t-s)}{R^2} +
(2+2\alpha)\int_0^t\gamma(z)dz\right\}}\gamma(s)\int_s^t[\nax u (\cdot,z)]_{\alpha}dzds \\
+C\|\rho_0\|_{L^{\infty}(\Rr^d)}\int_0^t\exp{\left\{-\fr{2\e(t-s)}{R^2} + 2\int_s^t\gamma(z)dz + 2\alpha\int_0^t\gamma(z)dz\right\}}[\nax u (\cdot, s)]_{\alpha}ds\\
+ C\left [[\tau_0]_{\alpha} + \|\tau_0\|_{L^{\infty}(\Rr^d)}\int_0^t[\nax u(\cdot,s)]_{\alpha}dt \right]\exp{\left\{-\fr{2 \e t}{R^2} +
(2+2\alpha)\int_0^t\gamma(s)ds\right\}}
\ea
\la{wsuho}
\ee
with $\gamma(t) = \sup_{x}|\nax u(x,t)|$ and the constant $C>0$ depending only on the dimension of space $d$.
\end{prop}
\noindent{\bf Proof.} The $L^{\infty}$ estimate follows from (\ref{qinf}) and
\be
\|e^i(\cdot, t)\|_{L^{\infty}(dx)} \le e^{\int_0^t\gamma(s)ds}
\la{eifty}
\ee
which follows immediately from (\ref{eieq}) and the fact that the initial data are of unit norm. The H\"{o}lder seminorm estimate follows from the algebra inequality
$$
[\phi\psi]_{\alpha} \le [\phi]_{\alpha}\|\psi\|_{L^{\infty}} + \|\phi\|_{L^{\infty}}[\psi]_{\alpha},
$$
(\ref{qinf}) and (\ref{eifty}) above, (\ref{einaxu}), (\ref{qalpha}), and the behavior of the H\"{o}lder continuous function under composition with $A$ and (\ref{alip}). The $W^{1,p}$ estimate follows using the product rule. This concludes the proof of the proposition.

Before we conclude this section, we comment on the case of $d=2$. In this case the velocity gradient  matrix is
\be
\nax u(x,t) = \left (
\begin{array}{cc}
\lambda(x,t) & \mu(x,t) -\frac{\omega(x,t)}{2}\\
\mu(x,t) + \frac{\omega(x,t)}{2} & -\lambda(x,t)
\end{array}
\right)
\la{nabu}
\ee
where
\be
\lambda(x,t) = \frac{1}{2}\left(\partial_1 u^1(x,t) -\partial_2u^2(x,t)\right),
\la{lam}
\ee
\be
\mu(x,t) = \frac{1}{2}\left (\partial_1 u^2(x,t) + \partial_2 u^1(x,t)\right) 
\la{mu}
\ee
and
\be
\omega(x,t) = \partial_1 u^2(x,t) -\partial_2u^1(x,t).
\la{ome}
\ee
The symmetric part of the gradient (rate of strain matrix) is
\be
S(x,t) = \left (
\begin{array}{cc}
\lambda(x,t) & \mu(x,t) \\
\mu (x,t)   & -\lambda(x,t)
\end{array}
\right )
\la{s}
\ee
Differentiating (\ref{eieq}) we can write
\be
D_t^2 e^{i} = \left((\nax u)^2 + D_t(\nax u)\right)e^{i}
\la{seqei}
\ee
This equation holds in any number of dimensions, but
because of the Hamilton-Cayley theorem, in two dimensions
we have
\be
(\nax u(x,t))^2 = \delta(x,t)\mathbb I
\la{hc}
\ee
with 
\be
\delta(x,t) = \lambda^2(x,t) + \mu^2(x,t) - \fr{1}{4}\omega^2(x,t) = -\det(\nax u)(x,t).\la{del}
\ee
Therefore the equations (\ref{seqei}) become
\be
D_t^2e^i = \left(\delta {\mathbb I} + D_t(\nax u)\right)e^i
\la{seqeid}
\ee
Similarly, we have
\be
D_s^2 Q = -Q\left(\delta{\mathbb I} + D_s(\nax u)\right) 
\la{seqeq}
\ee
We can easily integrate the equations (\ref{seqeid}, \ref{seqeq}) when the velocity gradient is constant in space and time (so that $D_t(\nax u)=0)$. In this case $\delta$ is a constant and its sign dictates the size of solutions $e^{i}$and $Q$. Consequently, when $\delta \le 0$ or when $\delta>0$ but $\sqrt{\delta}\le \frac{\e}{R^2}$, then 
$\tau$ is bounded, and when $\delta >(\fr{\e}{R^2})^2$, then $\tau$ may grow exponentially in time. Integrating (\ref{seqeid}, \ref{seqeq}) when $D_t(\nax u) = 0$ we see that, if $\delta \le 0$ the functions $e^{i}$ and $Q$ are bounded in time and when $\delta>0$, they grow like $\exp{\sqrt{\delta t}}$ in time.  If the initial distribution $f_0$ is radially symmetric in $m$ then the matrix $\sigma_0(a)$ is a multiple of the identity matrix (a scalar matrix, or an isotropic stress),
$$
\sigma_0(A(x,t)) = c\rho_0(A(x,t))\mathbb I.
$$
If $\delta$ is a positive constant, integrating (\ref{seqeid}) we have
$$
e^1(x,t) = \left(
\ba 
\cosh(t\sqrt{\delta})\\
\sinh(t\sqrt{\delta})
\ea
\right)
\quad\quad
e^2(x,t) = \left (
\ba
\sinh(t\sqrt{\delta})\\
\cosh(t\sqrt{\delta})
\ea
\right)
$$
In this case 
$$
\ba
e^{11}(x,t)= e^{22}(x,t) =c\rho_0(A(x,t))\cosh(2t\sqrt{\delta})\\ 
e^{12}(x,t) = c\rho_0(A(x,t))\sinh(2t\sqrt{\delta})
\ea
$$
The case of a uniform gradient is used in rheology, but if the velocity gradient is square-integrable in space, then $\delta$ has to have average zero in space because 
\be
(\nax u)^2_{ij} =\pa_k(u^{i}\pa_ju^k).
\la{avena}
\ee 
Therefore, for moderately varying $\nax u$ we can expect $\tau$ to be bounded.

\section{Coupling to the Stokes system}
The Stokes system is
\be
-\Dx u +\nax p = k\dx\sigma
\la{stks}
\ee
coupled with
\be
\dx u = 0
\la{dxu}
\ee
The coefficient $k$ has units of inverse time, and represents the ratio between an energy (kT) per unit mass and kinematic viscosity. The fluid density is taken to be one. We discuss the situation in which the fluid occupies all $\Rr^d$ and is at rest at infinity. The velocity and pressure decay in space and are obtained then from $\sigma $ by classical singular integrals. We can equally consider the case in which the velocity and pressure are periodic in space. Let us denote by
\be
R_j = \Lambda^{-1} \pa_j
\la{riesz}
\ee
for $j=1, \dots d$, in ${\Rr^d}$ the Riesz operators where
\be
\Lambda = (-\Delta)^{\fr{1}{2}}
\la{zy}
\ee
is the Zygmund operator. By modifying the pressure, we may write
\be
-\Dx u + \nax{\widetilde{p}} = k\dx(\tau)
\la{stws}
\ee
and using (\ref{dxu}) to solve for $\widetilde{p}$
\be
\widetilde{p} = -k R_mR_n(\tau^{mn})
\la{wp}
\ee
we deduce
\be 
u^{i} = k\Lambda^{-1}\left(R_l(\tau^{il}) + R_iR_mR_n(\tau^{mn})\right)
\la{ui}
\ee
and
\be
\pa_j u^i = kR_j\left(R_l(\tau^{il}) + R_iR_mR_n(\tau^{mn})\right)
\la{naws}
\ee
which we abbreviate as
\be
\nax u = k{\mathcal R}\tau
\la{nas}
\ee
The system formed by (\ref{qa}, \ref{ala})
together with (\ref{ui}) is closed. 

\subsection{Energetics}
The system formed by (\ref{lfp}) coupled to (\ref{stks}) has a Liapunov functional.
We take  (\ref{lfp}), multiply by $k(\log f + U)$,  integrate, take (\ref{stks}) multiply by $u$, integrate, and add the two: we deduce
\be
\ba
k\frac{d}{dt}\left[\int_{\Rr^d}\int_{\Rr^d}f\log fdm dx +\frac{1}{2}\int_{\Rr^d}Tr[\sigma(x,t)]dx \right]\\
= -\int_{\Rr^d}|\nax u(x,t)|^2dx -k\e\int_{\Rr^d}\int_{\Rr^d}f\left| \nam(\log f+U)\right|^2dmdx
\ea
\la{ens}
\ee
This is a general property of the class of Fokker-Planck equations we are interested in, but in  the case of the Oldroyd B system we are considering here, we can obtain the energetics directly from the equation (\ref{sigeq}) by taking the trace,  adding to the Stokes equation (\ref{stks}) multiplied by $\fr{2}{k} u$ and integrating:
\be
\pa_t\int_{\Rr^d} Tr\;\sigma dx + \fr{2}{k}\int_{\Rr^d}|\nax u|^2dx = -
\fr {2\e }{R^2}\int_{\Rr^d} Tr\;\sigma dx + \fr{2\e }{R^2}\int_{\Rr^d}\rho_0dx
\la{trsigeq}
\ee
Using the fact that $Tr\sigma$ together with $\rho$ bound the entries in $\tau$, we derive
\be
\|\tau(\cdot, t)\|_{L^1(\Rr^d)} \le C[\|\tau_0\|_{L^1(\Rr^d)} + \|\rho_0\|_{L^1(\Rr^d)}]
\la{wslone}
\ee
with a constant $C$ depending only on $d$.

\subsection{Gradient bounds}
We are going to investigate bounds on 
the time integral of the maximum gradient of velocity,
$$
\int_0^{T}\|\nax u(\cdot,t)\|_{L^{\infty}(dx)}dt
$$
In order to do so we recall a classical fact (\cite{stein}, Chapter 3, section 5), namely that operators given at Fourier level by multiplication by functions that are homogeneous of degree zero and are smooth on the unit sphere can be represented as sums of multiples of the identity and singular integral operators of a classical type. 
\begin{prop} {\la{rprop}}There exists a constant $C$ depending only on $d$ and $\alpha$ such that, for any $\ws$,  
\be
\|{\mathcal R}\ws\|_{L^{\infty}(dx)}\le C\|\ws\|_{L^{\infty}(dx)}\left \{1+ \log\left[1 + \fr{\|\ws\|_{L^1(\Rr^d)}^{\fr{\alpha}{d+\alpha}}[\ws]_{\alpha}^{\fr{d}{d+\alpha}}}{\|\ws\|_{L^{\infty}(dx)}}\right ]\right\}
\la{linfty}
\ee
where
\be
[\ws]_{\alpha} = \max_{mn}\sup_{x\neq y}\fr{|\ws^{mn}(x)-\ws^{mn}(y)|}{|x-y|^{\alpha}}
\la{semi}
\ee
and $0<\alpha<1$. There exists a constant $C$ depending only on $d$ and $\alpha$ such that
\be
[\mathcal R \ws]_{\alpha} \le C[\ws]_{\alpha}\la{wsalpha}
\ee
\end{prop} 
\noindent{\bf Proof.} In view of the fact that ${\mathcal R}$ is made from operators of the form $R_jR_k$ and $R_jR_kR_l$ we have the representation  
$$
({\mathcal R}\ws)_{ij}(x) = C_{ijmn}\ws^{mn}(x) + P.V.\int_{\Rr^d}\fr{k_{ijmn}(\hat{y})}{|y|^d}\ws^{mn}(x-y)dy
$$
where $C_{ijmn}$ are constants and each $k_{ijmn}$ is a smooth function (actually a harmonic polynomial of degree less or equal to three) of $\hat{y} = y/|y|$, with mean zero on the unit sphere ${\mathbb S}^{d-1}$. We need to prove the inequality for each entry, and we reduce the proof to the scalar case. The proof then follows along very classical lines. We break the integral
$$
K\ws (x) = P.V. \int_{\Rr^d}\frac{k(\hat y)}{|y|^d}\ws(x-y)dy
$$
into three pieces. We choose two numbers $m<M$ and write
$$
I_1(x) = P.V. \int_{|y|\le m}\frac{k(\hat y)}{|y|^d}\ws(x-y)dy
$$
$$
I_2(x) = \int_{m\le |y|\le M}\frac{k(\hat y)}{|y|^d}\ws(x-y)dy
$$
and 
$$
I_3(x) = \int_{|y|\ge M}\frac{k(\hat y)}{|y|^d}\ws(x-y)dy.
$$
We have, with the cancellation property of $k(\hat{y})$,
$$
|I_1(x)| \le Cm^{\alpha}[\ws]_{\alpha},
$$
then clearly
$$
|I_2(x)| \le C\log\left (\fr{M}{m}\right )\|\ws\|_{L^{\infty}}
$$
and 
$$
|I_3(x)|\le CM^{-d}\|\ws\|_{L^1(dx)}
$$
Now we choose the two length scales,
$$
M= \left [\fr{\|\ws\|_{L^1(\Rr^d)}}{\|\ws\|_{L^{\infty}(dx)}}\right]^{\fr{1}{d}}
$$
and
$$
m = \left [\fr{\|\ws\|_{L^{\infty}(dx)}}{[\ws]_{\alpha}}\right]^{\fr{1}{\alpha}}
$$
If $M\ge m$ we use both in the inequalities above and we obtain
$$
|K\ws(x)| \le C\|\ws\|_{L^{\infty}(dx)}\left[2 + \log\left(\fr{M}{m}\right)\right]
$$
If $M\le m$ then we forgo the division above and split the integrals into only two pieces, up to $M$ and from $M$. We obtain then
$$
|K\ws(x)| \le 2C\|\ws\|_{L^{\infty}(dx)}
$$
This ends the proof of (\ref{linfty}). The inequality (\ref{wsalpha}) is a classical inequality for singular integrals of the type above. This ends the proof of the proposition.
\section{Small data}
The approach we choose to control $[\tau]_{\alpha}$ is via the explicit formula (\ref{tau}). The equation obeyed by $\tau$ can be derived directly from (\ref{sigeq}):
\be
D_t\tau = -\fr{2\e}{R^2}\tau + (\nax u)\tau + \tau(\nax u)^T + 2\rho S
\la{taueq}
\ee
Because $2S = k{\mathcal R}\tau +k({\mathcal R}\tau)^T$, depends linearly on $\tau$, it is clear that strong enough damping will prevail, sending the solution to zero. It is clear however that
\be
\fr{\e}{kR^2}\ge C\|\rho_0\|_{L^{\infty}(\Rr^d)}
\la{ante}
\ee
will have to be part of the requirement for this stability result. The system of (\ref{lfp}) with (\ref{um}) and (\ref{sigmam}) coupled with (\ref{stks}) is equivalent to the system formed by (\ref{taueq}), the Stokes system (\ref{stws}) and (\ref{dtrho}). We use Proposition (\ref{pasprop}). 
We first prove a good local existence result. Let us denote
\be
M_1 = \|\rho_0\|_{L^1(\Rr^d)} + \|\tau_0\|_{L^1(\Rr^d)}
\la{mone}
\ee
\be
M_{\infty} = \|\rho_0\|_{L^{\infty}(\Rr^d)} + \|\tau_0\|_{L^{\infty}(\Rr^d)}
\la{mifty}
\ee
and 
\be
M_{\alpha} = [\rho_0]_{\alpha} + [\tau_0]_{\alpha}
\la{malpha}
\ee
Let us introduce the Deborah number
\be
D:= \frac{k}{\kappa_0} = \frac{kR^2}{\e},
\la{varepsilon}
\ee
the nondimensional ratio of the particle time scale $\frac{R^2}{\e}$
to the viscous fluid's response time to the stress added by the particles, 
$k^{-1}$.  We denoted half the damping rate by
\be
\kappa_0 = \frac{\e}{R^2}.
\la{kz}
\ee

\begin{thm}\la{lw} Let $(\tau_0, \rho_0)\in L^{1}(\Rr^d)^2\cap C^{\alpha}(\Rr^d)^2$.
There exist  constants $\varepsilon>0$, $\Gamma\ge 2$  and
there exists a time $T_0>0$ and a weak solution of (\ref{dtrho}, \ref{ui}, \ref{taueq})
$(\rho, \tau)\in C([0,T_0], W^{-1,p}(\Rr^d)^2)\cap L^{\infty}([0,T_0]; (C^{\alpha}(\Rr^d)\cap L^1(\Rr^d))^2)$, $1<p<\infty$, satisfying the equations in weak sense,  and such that
\be
Dk_0T_0 M_{\infty}\left\{1 + \log\left(1 + M_{\alpha}^{\fr{d}{d+\alpha}}M_1^{\fr{\alpha}{d+\alpha}}M_{\infty}^{-1}\right)\right\}\ge \varepsilon\la{tlow}
\ee
and
\be
\|\rho (\cdot,t)\|_{L^{1}(\Rr^d)} + \|\tau(\cdot, t)\|_{L^{1}(\Rr^d)} \le \Gamma M_{1},
\la{doubone}
\ee

\be
\|\rho (\cdot,t)\|_{L^{\infty}(\Rr^d)} + \|\tau(\cdot, t)\|_{L^{\infty}(\Rr^d)} \le \Gamma M_{\infty}
\la{doubifty}
\ee
and 
\be
[\rho (\cdot,t)]_{\alpha} + [\tau(\cdot, t)]_{\alpha} \le \Gamma M_{\alpha}
\la{doubalpha}
\ee
hold on $[0,T_0]$
\end{thm}
\noindent{\bf Proof.} We consider the system (\ref{taueq}, \ref{ui}, \ref{dtrho}) with initial data that have been regularized by convolution with a standard mollifier $\phi_{\delta}$. We obtain uniform bounds for quantities of interest and pass to the limit $\delta\to 0$, removing the mollifier. In order to simplify the exposition, we will denote the solutions $\tau^{(\delta)}, \rho^{(\delta)}$ by $\tau, \rho$ as if they did not depend on $\delta$. The bounds and constants below will indeed be $\delta$-independent.  We use the following notations
\be
{\og}(t) := (D\kappa_0)^{-1}\sup_{0\le s\le t}\|\nax u(\cdot, s)\|_{L^{\infty}(\Rr^d)},
\la{og}
\ee
\be
b(t):= (D\kappa_0)^{-1}\sup_{0\le s\le t}[\nax u(\cdot, s)]_{\alpha},\la{bt}
\ee
\be
n(t) : = \sup_{0\le s\le t}[\|\tau(\cdot, s)\|_{L^{\infty}(\Rr^d)} +\|\rho(\cdot, s)\|_{L^{\infty}(\Rr^d)}],\la{nt}
\ee
and 
\be
m(t): = \sup_{0\le s\le t}\{ [\tau(\cdot,s)]_{\alpha} + [\rho(\cdot, s)]_{\alpha}\}\la{mt}
\ee
We prove first uniform bounds, then discuss their use. We use (\ref{wslone})
\be
\sup_{0\le s\le t}[\|\tau (\cdot, s)\|_{L^1(\Rr^d)} + \|\rho (\cdot,s)\|_{L^1(\Rr^d)}\le C M_1\la{lonew}
\ee
which is based on the fundamental dissipative property of the system. We note that the function $(x,a,b)\mapsto x\log\left(1 +\fr{ab}{x}\right)$
is increasing in each of its arguments $x>0,\, a>0,\, b>0$. Therefore, we may use upper bounds in the right-hand side of (\ref{linfty}). 
Using (\ref{linfty}) we deduce
\be
\og(t) \le Cn(t)\left\{1 +\log\left (1 + \fr{m(t)^{\fr{d}{d+\alpha}}M_1^{\fr{\alpha}{d+\alpha}}}{n(t)}\right)\right\}
\la{ogbound}
\ee
and using (\ref{wsalpha})
\be
b(t) \le C m(t).
\la{btmt}
\ee
We proceed to estimate in (\ref{wsufty}):
\be
n(t)\le CM_{\infty}[1+D\kappa_0 t\og(t)]e^{CD\kappa_0 t\og(t)}\la{ntog}
\ee
and in (\ref{wsuho})
\be
m(t) \le CM_{\alpha}[1 + D\kappa_0 t\og(t)]e^{CD\kappa_0 t\og(t)} + C(M_{\infty}D\kappa_0 t)[1+ D\kappa_0 t\og(t)]e^{CD\kappa_0 t\og(t)}m(t)\la{mtog}
\ee
The term $D\kappa_0 t\og(t)$ can be bound from (\ref{ogbound})
\be
D\kappa_0 t\og(t) \le (CD\kappa_0 t)n(t)\left\{1 +\log\left(1 +\fr{m(t)^{\fr{d}{d+\alpha}}M_1^{\fr{\alpha}{d+\alpha}}}{n(t)}\right)\right\}\la{dkapog}
\ee
The structure of the inequalities is the following: if we denote the group
\be
x(t): = D\kappa_0 t\og(t), \la{xt}
\ee 
we have from (\ref{ntog})
\be 
n(t) \le CM_{\infty}(1+ x(t))e^{Cx(t)}\la{nx}
\ee
and therefore, from (\ref{dkapog}) we deduce
\be
\ba
x(t)\le (CD\kappa_0t)M_{\infty}(1+x(t))e^{Cx(t)}(1+\log (1+\fr{1}{C(1+x(t))e^{Cx(t)}}))\\\times\left\{ 1+ \log\left(1 + m(t)^{\fr{d}{d+\alpha}}M_1^{\fr{\alpha}{d+\alpha}}M_{\infty}^{-1}\right)\right\}
\ea
\la{temp}
\ee
so
\be
F(x(t))\le C(M_{\infty}D\kappa_0t)\left\{ 1+ \log\left(1 + m(t)^{\fr{d}{d+\alpha}}M_1^{\fr{\alpha}{d+\alpha}}M_{\infty}^{-1}\right)\right\}
\la{xbound}
\ee
where
\be
F(x) = xe^{-Cx}(1+x)^{-1}\left(1+\log (1+ \fr{1}{C(1+x)e^{Cx}})\right)^{-1}.
\la{Fx}
\ee
Let
\be
a(t): = M_{\infty}D\kappa_0 t
\la{at}
\ee
Now (\ref{mtog}) can be written as
\be
m(t)\le C(1+(x(t))e^{Cx(t)}[M_{\alpha} + a(t) m(t)]\la{mat}
\ee
and
\be
F(x)\le Ca(t)\left\{ 1+ \log\left(1 + m(t)^{\fr{d}{d+\alpha}}M_1^{\fr{\alpha}{d+\alpha}}M_{\infty}^{-1}\right)\right\}
\la{Fat}
\ee
We fix $C_0$ a large enough absolute constant, larger than the $C$ we encountered so far. Without loss of generality we may assume
$1+x\le e^{C_0x}$ and 
$$
1+\log (1+ \fr{1}{C(1+x)e^{Cx}})\le 1+\log(1+C^{-1})\le C_0
$$
so that
$$
F(x)\ge xe^{-2C_0x}\fr{1}{C_0}
$$
and therefore, as long as $x(t)\le 1$ we have from (\ref{Fat})
\be
x(t) \le C_0^2e^{2C_0}a(t)\left\{ 1+ \log\left(1 + m(t)^{\fr{d}{d+\alpha}}M_1^{\fr{\alpha}{d+\alpha}}M_{\infty}^{-1}\right)\right\}
\la{xm}
\ee
and, from (\ref{mat})
\be
\fr{m(t)}{M_{\alpha}} \le C_0e^{2C_0x(t)}\left [1 + a(t)\fr{m(t)}{M_{\alpha}}\right ].
\la{mm}
\ee
Using $x(t)\le 1$ and requiring
\be
a(t) \le \fr{1}{2}C_0^{-1}e^{-2C_0}\la{reqaone}
\ee
we deduce from (\ref{mm}) that, as long as $x(t)\le 1$, it follows that
\be
\fr{m(t)}{M_{\alpha}} \le 2C_0e^{2C_0}.
\la{mM}
\ee
Returning to $x(t)$ we have that, as long as $x(t)\le 1$, it follows from (\ref{xm}) and (\ref{mM}) that
\be
x(t) \le C_0^2e^{2C_0}a(t)\left\{1+ \log\left(1 +(2C_0e^{2C_0}M_{\alpha})^{\fr{d}{d+\alpha}}M_1^{\fr{\alpha}{d+\alpha}}M_{\infty}^{-1}\right)\right\}
\la{xtc}
\ee
We require therefore also
\be
a(t)\le \fr{1}{2}\left\{ C_0^2e^{2C_0}\left\{1+ \log\left(1 +(2C_0e^{2C_0}M_{\alpha})^{\fr{d}{d+\alpha}}M_1^{\fr{\alpha}{d+\alpha}}M_{\infty}^{-1}\right)\right\}\right\}^{-1}
\la{af}
\ee
and thus, as long as $x(t)\le 1$, it follows that actually $x$ is  bounded by $x(t)\le \fr{1}{2}$.
The initial datum for $x(t)$ is zero, and $a(t)$ is explicitly proportional to $t$ and starts from zero. Therefore, because of continuity, if the time $t$ is taken small enough so that $a(t)$ satisfies the requirements (\ref{reqaone}, \ref{af}), it follows that $x(t)<\fr{3}{4}$ on that interval, because we can reason by contradiction, and the first time it would reach $\fr{3}{4}$ it would have to be not larger that $\fr{1}{2}$, which is absurd.
Consequently, it follows from (\ref{mM}) that $m(t)\le 2C_0e^{2C_0}M_{\alpha}$ and from (\ref{nx}) that $n(t)\le 2C_0e^{2C_0}M_{\infty}$ on the time interval, providing the required short time uniform bound.
The rest of the proof of the theorem is based on the uniform bound. It is well-known that regular solutions exist as long as the $L^{\infty}$ bounds persist (see \cite{chemm}, \cite{leim}, also next section). Therefeore classical solutions with regularized initial data exist for this interval of time. Indeed, by regularizing the initial data we do not hamper the bounds $M_1, M_{\infty}, M_{\alpha}$. The most straightforward proof of local existence for smooth data yields a very short time of existence, but the a priori bounds permit the solution to be extended as long as (\ref{reqaone}) and (\ref{af}) hold. Then we remove the regularization of the initial data, and pass to limit in the equation in distribution sense.  We conclude that the limit obeys the equation, because the quadratic nonlinearity is strongly continuous in $L^2$ (for the right-hand side) and in $W^{-1,2}$ for the advective derivative. Further details are left for the interested reader.

We will pursue now the stability issue. Let us start by assuming that
\be
\fr{1}{\kappa_0}\|\nax u(\cdot, t)\|_{L^{\infty}(\Rr^d)} \le Ge^{-\kappa_0 t}
\la{ass}
\ee
on the time interval $t\in [0,T]$. 
We require that
\be
{G}\le \fr{1}{4}.
\la{kappazero}
\ee
Similarly, we assume
\be
\fr{1}{\kappa_0}[\nax u(\cdot, t)]_{\alpha} \le He^{-\kappa_0 t}
\la{assu}
\ee
on the same time interval.
We will show that we can choose $G$, $H$ so that these assumptions are invariant in time: once initiated they cannot break down. Before we embark on the calculation, we should point out the simple strategy: we obtain from (\ref{ass}), (\ref{assu}) a priori bounds on $\tau$. We then deduce a posteriori bounds for $\nax u$ that are strictly better than the assumptions
(\ref{ass}, \ref{assu}). Because we work in a class in which the equation has unique solutions, and if the initial data for $\tau$ are small, we deduce that, if initiated, the  bounds continue indefinitely.

We start by estimating the terms in (\ref{wsuho}) using (\ref{ass}) and (\ref{assu}). 
We will use
\be
\int_0^t \gamma(z) dz \le  G
\la{G}
\ee
which follows from (\ref{ass}).
We have four terms, and in all of them, as well as in (\ref{wsufty})
we bound, 
\be
\ba
\int_0^t\exp{\{-2\kappa_0(t-s) + 2\int_s^t\gamma(z)\}}\gamma(s) ds\\ \le
e^{2G}\exp{\{-2\kappa_0 t\}}\int_0^{t}\exp{\{2\kappa_0s\}}\kappa_0G\exp{\{-\kappa_0 s\}}ds \le {G}e^{-\kappa_0t}e^{2G}\\
\le \sqrt{e}Ge^{-\kappa_0 t}
\ea
\la{intbasis}
\ee
using (\ref{G}) and the fact that $4G\le 1$ which is (\ref{kappazero}).
From this calculation and (\ref{wsufty}) it follows that
\be
\|\tau (\cdot,t)\|_{L^{\infty}(\Rr^d)}\le Ce^{-\kappa_0 t}\left[\|\rho_0\|_{L^{\infty}(\Rr^d)} +\|\tau_0\|_{L^{\infty}(\Rr^d)}\right] = Ce^{-\kappa_0 t}M_{\infty}
\la{tauassfty}
\ee
with $C$ depending on $d$ alone. Using  (\ref{qinf}) and (\ref{eifty}) in (\ref{tau}) and integrating we obtain in the same way
\be
\|\tau (\cdot,t)\|_{L^{1}(\Rr^d)}\le Ce^{-\kappa_0 t}\left[\|\rho_0\|_{L^{1}(\Rr^d)} +\|\tau_0\|_{L^{1}(\Rr^d)}\right] = Ce^{-\kappa_0t}M_{1}
\la{tauassone}
\ee 
We write the inequality (\ref{wsuho}) as
\be
[\tau(\cdot, t)]_{\alpha} \le I+ II + III + IV
\la{ontofour}
\ee
and we start with
\be
I \le C[\rho_0]_{\alpha}e^{\alpha G}\left(\fr{1}{2}\right) e^{-\kappa_0 t}= C[\rho_0]_{\alpha}e^{-\kappa_0 t}.
\la{I}
\ee
in which we used (\ref{G}) and (\ref{intbasis})  above.
This illustrates our use of the generic constant $C$: we used $\alpha G\le \fr{1}{2}$.
For $II$ we obtain
\be
\ba
II \le C\|\rho_0\|_{L^{\infty}(\Rr^d)}e^{2\alpha G}e^{-\kappa_0 t}\int_0^t[\nax u(\cdot, s)]_{\alpha}ds\\
\le C\|\rho_0\|_{L^{\infty}(\Rr^d)}e^{-\kappa_0 t}\int_0^t[\nax u(\cdot, s)]_{\alpha}ds
\ea
\la{II}
\ee
We used ${2\alpha G}\le {1}$ and (\ref{intbasis}). For $III$ we need
to use the assumption (\ref{assu}) and instead of (\ref{intbasis}) we use
a similar bound with the explicit $\kappa_0G\exp{(-\kappa_0 t)}$ replacing $\gamma (t)$ and $\kappa_0 H\exp{(-\kappa_0 t)}$ replacing $[\nax u]_{\alpha}(t)$. We obtain:
\be
III \le C\|\rho_0\|_{L^{\infty}(\Rr^d)}e^{(2 + 2\alpha)G}He^{-\kappa_0 t} = C\|\rho_0\|_{L^{\infty}(\Rr^d)}He^{-\kappa_0 t} 
\la{III}
\ee
Finally,
\be
IV \le C\left[[\tau_0]_{\alpha} +\|\tau_0\|_{L^{\infty}(\Rr^d)}\int_0^t[\nax u(\cdot,s)]_{\alpha}ds\right]e^{-2\kappa_0 t}
\la{IV}
\ee
where we use again (\ref{G}) and $(2+2\alpha)G \le 1$. 
Collecting terms we have proved that  (\ref{ass}, \ref{kappazero})  
and (\ref{assu}) imply that
\be
[\tau (\cdot, t)]_{\alpha}\le Ce^{-\kappa_0 t}\left\{ [\rho_0]_{\alpha} + [\tau_0]_{\alpha} + H\|\rho_0\|_{L^{\infty}(\Rr^d)} +
\left [\|\rho_0\|_{L^{\infty}(\Rr^d)}+ \|\tau_0\|_{L^{\infty}(\Rr^d)}\right] \int_0^t[\nax u(\cdot,s)]_{\alpha}ds\right\}
\la{taualpha}
\ee
Now from (\ref{nas}), (\ref{wsalpha}) and  a Gronwall-like inequality we derive
\be
[\tau(\cdot, t)]_{\alpha} \le C\left[[\rho_0]_{\alpha} + [\tau_0]_{\alpha} + H\|\rho_0\|_{L^{\infty}(\Rr^d)}\right]\exp{\left \{-t\kappa_0 +
\fr{Ck[\|\rho_0\|_{L^{\infty}(\Rr^d)} + \|\tau_0\|_{L^{\infty}(\Rr^d)}]}{\kappa_0}\right\}}
\la{wstho}
\ee
Indeed, (\ref{taualpha}) together with (\ref{nas}, \ref{wsalpha}) imply an inequality of the type
$$
y(t) \le Ce^{-t\kappa_0}\left [a +b\int_0^ty(s)ds\right] = F(t)
$$
for the  positive quantity $y(t) = [\tau(\cdot, t)]_{\alpha}$, with 
$$
a=[\rho_0]_{\alpha} + [\tau_0]_{\alpha} + H\|\rho_0\|_{L^{\infty}(\Rr^d)}
$$
and  
$$
b=k[\|\rho_0\|_{L^{\infty}(\Rr^d)} + \|\tau_0\|_{L^{\infty}(\Rr^d)}].
$$
Differentiating the right hand side $F(t)$ and using the fact that $y(t)\le F(t)$
we deduce
$$
\fr{dF}{dt}\le -\kappa_0 F + Cbe^{-t\kappa_0} F
$$
dividing by $F$ and integrating we obtain
$$
\log\left(\fr{F(t)}{Ca}\right) \le -\kappa_0 t + \fr{Cb}{\kappa_0}
$$
and substituting back, we deduce (\ref{wstho}). Let us consider now the right-hand side of (\ref{linfty}) with $\ws = \tau$. 
We recall that (\ref{tauassfty}) is 
\be
\|\tau(\cdot, t)\|_{L^{\infty}(\Rr^d)} \le Ce^{-\kappa_0 t}M_{\infty},
\la{taufty}
\ee
(\ref{tauassone}) is
\be
\|\tau(\cdot, t)\|_{L^{1}(\Rr^d)} \le Ce^{-\kappa_0 t}M_{1},
\la{tauone}
\ee
and (\ref{wstho}) becomes
\be
[\tau (\cdot, t)]_{\alpha} \le C[M_{\alpha} +H\|\rho_0\|_{L^{\infty}(\Rr^d)}]e^{-\kappa_0 t}\exp{\{CDM_{\infty}\}}\la{tautalpha}
\ee
Assembling the right hand side of (\ref{linfty}) we deduce 
\be
\|{\mathcal R}\tau\|_{L^{\infty}(\Rr^d)}\le CM_{\infty}e^{-\kappa_0 t}\left\{ 1+ \log\left[1 +\fr{M_1^{\fr{\alpha}{d+\alpha}}[M_{\alpha} +HM_{\infty}]^{\fr{d}{d+\alpha}}}{M_{\infty}}\right] + DM_{\infty}\right\} 
\la{inter}
\ee
Note that the term $e^{-\kappa_0 t}$ cancelled in the expression for the argument of the logarithm.  Now we use (\ref{nas}) to obtain a posteriori estimates. From (\ref{inter}) we obtain 
\be
\ba
\fr{1}{\kappa_0}\|(\nax u(\cdot, t)\|_{L^{\infty}(\Rr^d)}\\ \le CDM_{\infty}e^{-\kappa_0 t}\left\{1 + \log\left [1 + \fr{M_1^{\fr{\alpha}{d+\alpha}}[M_{\alpha} + HM_{\infty}]^{\fr{d}{d+\alpha}}}{M_{\infty}}\right] + DM_{\infty}\right\}
\ea
\la{naxuf}
\ee
Also, from (\ref{tautalpha}) and (\ref{nas}) we have
\be
\fr{1}{\kappa_0}[(\nax u)(\cdot, t)]_{\alpha} \le CD[M_{\alpha} + HM_{\infty}]\exp{\{-\kappa_0t +CDM_{\infty}\}}
\la{naxalphaf}
\ee
Let us choose $G=1/4$ and
\be
H = \fr{M_{\alpha}}{M_{\infty}}\la{H},
\ee
and so (\ref{naxalphaf}) becomes
\be
\fr{1}{\kappa_0}[(\nax u)(\cdot, t)]_{\alpha} \le CDM_{\alpha}\exp{\{-\kappa_0t +CDM_{\infty}\}}
\la{naxualphad}
\ee
The bound (\ref{tautalpha}) becomes
\be
[\tau (\cdot, t)]_{\alpha} \le CM_{\alpha}e^{-\kappa_0 t}\exp{\{CDM_{\infty}\}}\la{tautalphad}
\ee
The upper bound (\ref{naxuf}) implies the  bound
\be
\fr{1}{\kappa_0}\|(\nax u)(\cdot, t)\|_{L^{\infty}(\Rr^d)}\le B_1Ge^{-\kappa_0 t}
\la{naxunewt}
\ee
where
\be
B_1= B_0 + CD^2M_{\infty}^2\la{b1}
\ee
and $B_0= B_0(M_1,M_{\alpha}, M_{\infty}, D)$ is given by 
\be
B_0(M_1,M_{\alpha}, M_{\infty}, D) = CDM_{\infty}\left\{1 + \log\left[1+ M_{\infty}^{-1}M_{\alpha}^{\fr{d}{d+\alpha}}M_1^{\fr{\alpha}{d+\alpha}}\right ]\right\} 
\la{bzero}
\ee
with appropriate constants $C$ depending on $d$ only. 
The functions $B_0(M_1,M_{\alpha}, M_{\infty}, D)$ and $B_1(M_1, M_{\alpha}, M_{\infty}, D)$ are  continuous in their
arguments $M_1>0$, $M_{\alpha}>0$, $M_{\infty}>0$, $D>0$, and  hence locally bounded, and vanish when $M_{\infty}\to 0$ while the rest are held fixed, or when $D\to 0$ and the rest are fixed. 

The bound (\ref{naxualphad}) implies the  bound
\be
\frac{1}{\kappa_0}[(\nax u)(\cdot, t)]_{\alpha} \le B_2 He^{-\kappa_0 t} 
\la{naxualphat}
\ee
where the function $B_2= B_2(DM_{\infty})$ is given by
\be
B_2 = C DM_{\infty}\exp{\{C DM_{\infty}\}}\la{b2}
\ee
with an appropriate constant $C$. We have proved therefore

\begin{prop}\la{kine} Let $\tau(x,t)$ be the solution of (\ref{taueq}) 
given in (\ref{tau}) where $\rho(x,t)$ solves (\ref{dtrho}) and assume that $u$ obeys
(\ref{ass}), (\ref{assu})) on the time interval $[0,T]$.  Then $\tau$ and $u$ obey the bounds (\ref{taufty}, \ref{tauone}, \ref{tautalphad}, \ref{naxunewt}, \ref{naxualphat}) one the same time interval, where the constants $M_1, M_{\infty}, M_{\alpha}$ are given in (\ref{mone}, \ref{mifty}, \ref{malpha}), the constants $B_0$, $B_1, B_2$ are given in (\ref{b1}, \ref{bzero}, \ref{b2}) and the constants $C$ are independent of solutions, fixed, and depend only on $d$. 
\end{prop}
We combine now this proposition with the good local existence result, Theorem \ref{lw}. We start with initial data that are small.
To fix the notation, we denote by $C_1$ an upper bound for the constants $C$ appearing in Proposition \ref{kine} above.
We can arrange the constants so that, from
\be
[\tau(\cdot, t) ]_{\alpha} \le \Gamma M_{\alpha},
\la{gaalpha}
\ee
\be
\|\tau (\cdot, t)\|_{L^{1}(\Rr^d)} \le \Gamma M_{1},                   
\la{gaone}
\ee
and
\be
\|\tau (\cdot, t)\|_{L^{\infty}(\Rr^d)} \le \Gamma M_{\infty}
\la{gafty}
\ee
it follows that
\be
\fr{1}{\kappa_0}\|\nax u(\cdot, t)\|_{L^{\infty}(\Rr^d)}\le \fr{G}{8C_1}
\la{naxgsm}
\ee
and
\be
\fr{1}{\kappa_0}[\nax u (\cdot, t)]_{\alpha}\le \fr{H}{8C_1}
\la{naxalsm}
\ee
because of the constitutive relations between $u$ and $\tau$. Note that the size of $\tau$ alone determines that of $u$. Indeed, in view of (\ref{linfty}),
$$
\fr{1}{\kappa_0}\|\nax u(\cdot, t)\|_{L^{\infty}(\Rr^d)} \le CB_0(\Gamma M_1, \Gamma M_{\alpha}, \Gamma M_{\infty}, D)
$$
and  in view of (\ref{wsalpha})
$$
\fr{1}{\kappa_0}[\nax u(\cdot, t)]_{\alpha}\le CD\Gamma M_{\alpha}= CD\Gamma M_{\infty} H 
$$
So, the condition
\be
B_0(\Gamma M_1, \Gamma M_{\alpha}, \Gamma M_{\infty}, D) \le \fr{1}{32CC_1}
\la{smone}
\ee
implies (\ref{naxgsm}) and the condition
\be
D\Gamma M_{\infty} \le\fr{1}{8CC_1}\la{smtwo}
\ee
implies (\ref{naxalsm}). Then, we have automatically that (\ref{ass}) and (\ref{assu}) hold on the interval of time $[0, T_0]$. In fact,  if  $e^{\kappa_0 T_0}\le 2C_1$ we deduce that
\be
\fr{1}{\kappa_0}\|\nax u (\cdot, t)\|_{L^{\infty}(\Rr^d)} < \fr{1}{2}Ge^{-\kappa_0 t}
\la{halfass}
\ee
and
\be
\fr{1}{\kappa_0}[\nax u (\cdot, t)]_{\alpha} < \fr{1}{2} He^{-\kappa_0t}
\la{halfassu}
\ee
hold on $[0, T_0]$. On the other hand, as long as 
(\ref{ass}, \ref{assu}) hold, the decay estimates (\ref{taufty}, \ref{tauone}, \ref{tautalphad}) imply that
\be
\|\tau (\cdot, t)\|_{L^{\infty}(\Rr^d)} \le \fr{2}{3}M_{\infty},
\la{rifty}
\ee
\be
\|\tau (\cdot, t)\|_{L^{1}(\Rr^d)} \le \fr{2}{3}M_{1}
\la{rone}
\ee
and 
\be
[\tau (\cdot, t)]_{\alpha} \le \fr{2}{3}M_{\alpha}
\la{ralpha}
\ee
provided $e^{\kappa_0 t} \ge \fr{3}{2} C_1e^{C_1DM_{\infty}}$. The lower bound (\ref{tlow}) on the life-span of the solution  is
\be
\kappa_0 T_0 B_0(M_1, M_{\alpha}, M_{\infty}, D)\ge \varepsilon
\la{tlowb}
\ee
and it shows that if
\be
B_0(M_1, M_{\alpha}, M_{\infty}, D)\left[\log\left(\fr{3C_1}{2}\right ) + C_1DM_{\infty}\right] \le \varepsilon
\la{smthree}
\ee
then the solution exists for a period of time that is long enough so that $\tau$  is reduced from its original size (\ref{rifty}, \ref{rone}, \ref{ralpha}) at the end $T_0$ of the interval. Requiring that
\be
C_1DM_{\infty} \le \log\left(\fr{4}{3}\right)
\la{smfour}
\ee
allows to arrange for $T_0$ so that $\fr{3}{2}C_1e^{C_1DM_{\infty}}\le e^{\kappa_0T_0}\le 2C_1$. Therefore, if the initial data satisfy (\ref{smone}, \ref{smtwo}, \ref{smthree}, \ref{smfour}) then the assumptions of Proposition \ref{kine} are satisfied with $T=T_0$. We can repeat now argument taking the same time step. Therefore we find that the solution, if unique, exists for all time and obeys the bounds (\ref{gaalpha},\ref{gaone}, \ref{gafty}) for all time. But then, inspecting the bounds (\ref{wsufty}) and (\ref{wsuho}) in Proposition {\ref{pasprop}} we deduce that $\tau$ decays exponentially, and consequently so does $u$.

We note that the smallness conditions (\ref{smone}, \ref{smtwo}, \ref{smthree}, \ref{smfour}) all follow from a single condition
\be
B_0(M_1, M_{\alpha}, M_{\infty}, D) \le\varepsilon_1
\la{sm}
\ee
\begin{thm}\la{wsm} 
Let the initial data $\tau_0, \rho_0$ satisfy $(\tau_0, \rho_0) \in (L^{1}(\Rr^d))^2\cap (C^{\alpha}(\Rr^d))^2$, $0<\alpha<1$. Let $M_1, M_{\infty}$ and $M_{\alpha}$ defined in (\ref{mone}, \ref{mifty}, \ref{malpha}) denote the size of the initial data. There exists a constant $\varepsilon_1$ such that, if 
\be
DM_{\infty}\left\{1+ \log\left(1 + M_{\infty}^{-1}M_{\alpha}^{\fr{d}{d+\alpha}}M_1^{\fr{\alpha}{d+\alpha}}\right)\right\} \le \varepsilon_1
\la{smeps}
\ee
then there exists a unique global weak solution of $(\tau, \rho) \in
L^{\infty}([0,\infty), (L^1(\Rr^d)\cap C^{\alpha}(\Rr^d))^2) \cap C([0,\infty), (W^{-1,p}(\Rr^d))^2)$, $p<\infty$ of
(\ref{dtrho}), (\ref{taueq}), with $u$ given by (\ref{ui}). 
The meaning of weak solutions is in the sense of distributions, and the time derivatives of $\tau$ and $\rho$ belong to $W^{-1,p}(\Rr^d)$ for any $p<\infty$.
The gradient of velocity decays exponentially in time
\be
\fr{1}{\kappa_0}\|\nax u(\cdot,t)\|_{L^{\infty}(\Rr^d)} \le Ke^{-\kappa_0 t}
\la{naxdecay}
\ee
\be
\fr{1}{\kappa_0}[\nax u(\cdot,t)]_{\alpha} \le Ke^{-\kappa_0 t}
\la{naxalphadecay}
\ee
The norms $[\tau]_{\alpha}$ and $\|\tau\|_{L^{p}(\Rr^d)}$ decay exponentially:
\be
[\tau (\cdot,t)]_{\alpha} \le Ke^{-\kappa_0 t},\la{wsadecay}
\ee
\be
\|\tau (\cdot, t)\|_{L^{p}(\Rr^d)} \le K_p e^{-\kappa_0 t}
\la{wsinftydecay}
\ee
for $1\le p\le \infty$ with explicit constants $K$ that depend on $d$, $D$, $M_1$, $M_{\infty}$, $M_{\alpha}$, $p$  and $\alpha>0$ alone.
The norms of $\rho$ remain bounded.
\be
\|\rho(\cdot ,t)\|_{L^p(\Rr^d)} = \|\rho_0\|_{L^p(\Rr^d)}\la{rhob}
\ee
and
\be
[\rho (\cdot, t)]_{\alpha}\le K. \la{rhoalb}
\ee
\end{thm}
\begin{rmk} The condition (\ref{smeps}) is satisfied if the Deborah number is arbitrary, the initial $\tau_0$ is sufficiently small in $C^{\alpha}(\Rr^d)$ and the initial $\rho_0$ is small in $C^{\alpha}$. Nontrivial initial data with small $L^{\infty}$ norm and large spatial derivatives are allowed as well, because of the logarithmic dependence on $M_{\alpha}$. The conditions are also satisfied if the initial data are of order one but the Deborah number is small. If the initial data are smoother, the smoothness propagates.
\end{rmk} 
\begin{rmk} The proof of uniqueness, given below, uses Lagrangian transformations. Unlike the case of Euler equations, a proof in Eulerian coordinates seems difficult to obtain.
\end{rmk}
\noindent{\bf Proof.} The previous argument is complete, modulo the uniqueness of the solution. We present here the proof of uniqueness. We consider the equation 
\be
\fr{d}{dt} X = F[X]\la{lag}
\ee
with $X(a,0) = a$, where
$X= X(a,t)$ is viewed as an element of ${\mathcal{X}} := C([0,T], C^{1,\alpha}(\Rr^d)^d)$. We fix $T$.
The function $F$ depends on the whole path $X$ not only on the value of $X$ at some point and is obtained as follows. First we construct, 
using the initial data $\rho_0\in L^{1}(\Rr^d)\cap C^{\alpha}(\Rr^d)$, 
$\sigma_0(a) = \rho_0(a){\mathbb I} + \tau_0(a) \in (L^{1}(\Rr^d)\cap C^{\alpha}(\Rr^d))^{d^2}$, and using  (\ref{all}), the map $X\mapsto \sigma[X]$ given by 
\be
\ba
\sigma[X](a,t)\\
 = 2\kappa_0\rho_0(a)\int_0^t e^{-2\kappa_0(t-s)}q(a,t,s)q(a,t,s)^Tds + e^{-2\kappa_0 t}(\naa X(a,t))\sigma_0(a)(\naa X(a,t))^T
\ea
\la{qt}
\ee
where $q(a,t,s)= q[X](a,t,s) = (\naa X(a,t))(\naa X(a,s))^{-1}$ is given in (\ref{qa}). We take 
\be
\tau[X](a,t) = \sigma[X](a,t) - \rho_0(a){\mathbb I},\la{taux}
\ee
consider $A(x,t) = X^{-1}(x,t)$, compose $\tau[X]$ with it,
\be
\tau_X(x,t) =\tau[X](X^{-1}(x,t),t),
\la{tauxa}
\ee
 and solve the Stokes system, resulting in (\ref{ui})
\be
u_X^i(x,t) = k\Lambda^{-1}(R_l\tau_X^{il} + R_iR_mR_n(\tau_X^{mn})).\la{ux}  
\ee
We write symbolically
\be
u_X = k\Lambda^{-1}{\mathbb H}(\tau[X]\circ X^{-1}) = k\Lambda^{-1}{\mathbb H}(\tau_X)
\la{uxx}
\ee
where ${\mathbb H}$ stands for the combinations of Riesz transforms that appear in (\ref{ux}) 
\be
{\mathbb H}_{imn} = \delta_{im}R_n + R_iR_mR_n\la{hop}
\ee
and thus $k\Lambda^{-1}\mathbb H$ is the inverse of the Stokes system (\ref{stks}).  Finally, we compose with $X(a,t)$
\be
F^i[X](a,t) = u_X^i(X(a,t),t).\la{fx}
\ee
Thus, $F[X]$ is obtained via the succession of compositions
\be 
X\mapsto \tau[X]\mapsto \tau_X= \tau[X]\circ X^{-1}\mapsto u_X = k\Lambda^{-1}{\mathbb H}(\tau_X)\mapsto F= u_X\circ X
\la{chain}
\ee
The norm in ${\mathcal X}$ is
\be
\| X\|_{\mathcal X} := \sup_{0\le t\le T}\|X(\cdot,t)\|_{C^{1,\alpha}(\Rr^d)}
\la{normx}
\ee
We consider a fixed constant $M$ and the set
\be
{\mathcal D}:=\{ X\in {\mathcal X}\;\left |\right .\; X(0,a) = a, \; \fr{1}{2}\le \det\naa X(a,t) \le\fr{3}{2},\;\| X\|_{\mathcal X} \le M\}
\la {mathcald}
\ee
The initial data for the PDE serve as parameters in the definition of $F$. We 
wish to show that two solutions $X_1\in {\mathcal D}$ and $X_2\in {\mathcal D}$ of the equation (\ref{lag}) correspondingto the same $\rho_0, \tau_0$, are identical. In order to do so we establish
\be
\| (DF[X])Y\|_{\mathcal X} \le C \|Y\|_{\mathcal X}
\la{lip}
\ee
with a uniform constant $C$ that depends on $M$. We have to be careful to avoid taking derivatives of $\sigma_0$ and $\tau_0$. We start by noting that the map
$$
X \mapsto \sigma[X]
$$
is Fr\'{e}chet differentiable at $X\in {\mathcal D}$ as a map from the Banach space ${\mathcal X}$ to the Banach space $\Sigma = C(0,T; [C^{\alpha}(\Rr^d)\cap L^1(\Rr^d)]^{d^2})$ of time-continuous maps with values in the space of matrices with spatially H\"{o}lder continuous and integrable coefficients.
The derivative is a bounded linear map in ${\mathcal L}({\mathcal X}, \Sigma)$,$$
Y\mapsto (D\sigma[X])Y.
$$
The derivative has a complicated expression that depends on history but it can be easily obtained. The derivative of $q[X]$ is 
\be
((Dq[X])Y)(a,t,s) =(\pa_a Y(a,t))(\naa X(a,s))^{-1} - q(a,t,s)(\naa Y(a,s))(\naa X(a,s))^{-1}.
\la{Dq}
\ee 
Then the derivative of $\sigma[X]$ is
\be
\ba
((D\sigma[X])(Y))(a,t) =
2\kappa_0\rho_0(a)\int_0^te^{-2\kappa_0(t-s)}(Dq[X]Y)(a,t,s)q(a,t,s)^Tds \\
+ 2\kappa_0\rho_0(a)\int_0^te^{-2\kappa_0(t-s)}q(a,t,s)(Dq[X]Y)^T(a,t,s)ds \\
+e^{-2\kappa_0 t}\left [(\naa Y(a,t))\sigma_0(a)(\naa X(a,t))^T +
(\naa X(a,t))\sigma_0(a)(\naa Y(a,t))^T\right].
\ea
\la{dsig}
\ee  
It is clear that
\be
\|D\sigma[X]Y\|_{\Sigma} \le C\|Y\|_{\mathcal X}
\la{dsigmab}
\ee
where
\be
\|\tau\|_{\Sigma} = \sup_{0\le t\le T}\left\{\|\tau(\cdot,t)\|_{L^{\infty}(\Rr^d)} + \|\tau (\cdot, t)\|_{L^1(\Rr^d)} + [\tau (\cdot,t)]_{\alpha}\right\}.
\la{ny}
\ee
The map 
$$
X\mapsto \tau[X]
$$
is just a translation of $\sigma[X]$ by an $X$-independent amount, so it has the same derivative 
\be
(D\tau[X])Y = (D\sigma[X])Y.
\la{deq}
\ee
Now, the map 
$$
X(a,t)\mapsto \tau_X(x,t)
$$
is obtained from the map $\tau[X]$ by composition with $X^{-1}$. We cannot afford to take the derivative of $\tau_X$. However, $u_X$ is obtained from 
$\tau_X$ using a linear smoothing operator of degree minus one, so we can differentiate it. 
We take a path $X_{\varepsilon}(a,t) = X(a,t) + \varepsilon Y(a,t)$, differentiate with respect to $\varepsilon$ and then set $\varepsilon =0$:
\be
((Du_X)Y)(x,t) = k\Lambda^{-1}\mathbb H((D\tau[X]Y)\circ X^{-1}) - 
k\Lambda^{-1}{\mathbb H}((\nax \tau_X)(Y\circ{X^{-1}}))
\la{dux}
\ee
The first term  appears when we differentiate $\tau[X]\circ X^{-1}$ with respect to $X$, keeping $X^{-1}$ fixed. The second term arises when we differentiate $X^{-1}$ using the fact that
$$
\fr{d}{d\varepsilon}{X_{\varepsilon}^{-1}}_{\left | \right . \varepsilon = 0}(x,t) = - (\nax X^{-1})Y(X^{-1}(x,t),t).
$$
This fact is obtained by differentiating $\fr{d}{d\varepsilon}$
$$
X_{\varepsilon}^{-1}(X(a,t) +\varepsilon Y(a,t)) = a,
$$
setting $\varepsilon =0$, and reading at $a=X^{-1}(x,t)$. Then the derivative of $\tau[X]\circ X_{\varepsilon}^{-1}$  with respect to $\varepsilon$ is obtained using the chain rule
$$
\ba
\fr{d}{d\varepsilon}(\tau[X])(X_{\varepsilon}^{-1})_{\left |\right .\varepsilon =0} = ((\naa \tau[X])\circ X^{-1})(\fr{d}{d\varepsilon}X_{\varepsilon}^{-1})_{\left | \right . \varepsilon = 0}\\
= -((\naa \tau[X])\circ X^{-1})(\nax X^{-1})(Y\circ X^{-1}) = -(\nax \tau_X)(Y\circ X^{-1})
\ea
$$
The second term in (\ref{dux}) is strictly speaking a distribution, as it involves derivatives of $\tau_X$ which is only H\"{o}lder continuous. At this stage, we can view this as a formal calculation that will be justified in the end.
Finally, we need to compose back with $X(a,t)$.
\be
(DF(X)Y)(a,t) = ((Du_X)Y)(X(a,t),t) + (\nax u_X)(X(a,t),t)Y(a,t)
\la{df}
\ee
The second term in the expression above is {\em{unbounded}} as a linear operator from $C^{1,\alpha}$ to itself, simply because the coefficient $\nax u_X$ is not differentiable. But this term combined with the second problematic term in (\ref{dux}) produces a commutator that is better behaved, $K[X]Y$. Thus, we have
\be
(DF(X)Y)(a,t) = k\left\{\Lambda^{-1}{\mathbb H}((D\tau[X]Y)\circ A)\right\}(X(a,t), t)\la{dfx} + (K[X]Y)(a,t)
\ee
where
\be
(K[X]Y)(a,t) = (\nax u_X)(X(a,t),t)Y(a,t) - k(\Lambda^{-1}{\mathbb H}((\nax\tau_X)(Y\circ X^{-1})))(X(a,t),t)
\la{kx}
\ee
We would like to show that the map $Y\mapsto K[X]Y$ is a bounded linear operator in $C(0,T;(C^{1,\alpha})^d)$ with norm uniformly bounded for $X\in {\mathcal D}$. Becuse composition with $X^{-1}$ and composition with $X$ are both bounded linear operators $C^{1,\alpha}\to C^{1,\alpha}$, with norms controlled by $M$, the boundedness of $K[X]Y$ is equivalent to the boundedness of the map
$$
\phi \mapsto L[X]\phi 
$$
where
\be
\phi(x,t) = Y(X^{-1}(x,t),t)
\la{phix}
\ee
and
\be
(L[X]\phi)(x,t) = k(\nax \Lambda^{-1}{\mathbb H}\tau_X)(x,t)\phi(x,t) - k(\Lambda^{-1}{\mathbb H}((\nax \tau_X)\phi))(x,t).
 \la{lx}
\ee
It is important to specify the tensorial nature of this commutator, as not all
such expressions are better behaved than their individual terms. In our case,
$\tau_X$ is a fixed symmetric matrix in $\Sigma$, $k$ is a number, and the commutator is
\be
(L[X]\phi)^i = k \left[ \phi^p\pa_p\Lambda^{-1}{\mathbb H}_{imn}\tau_X^{mn} -
\Lambda^{-1}{\mathbb H}_{imn}(\phi_p\pa_p\tau_X^{mn})\right ]
\la{comcom}
\ee
We can write the commutator as
\be
L[X]\phi = k\left[\phi^p(\pa_p\Lambda^{-1}{\mathbb H}_{imn})(\tau_X^{mn})-(\pa_p\Lambda^{-1}{\mathbb H}_{imn})(\phi^p\tau_X^{mn})\right] -k\Lambda^{-1}{\mathbb H}_{imn}((\pa_p\phi^p)\tau_X^{mn})
\la{comsum}
\ee
For incompressible $X$  we could only consider divergence-free $\phi$ 
\be
\pa_p \phi^p = 0
\la{divphi}
\ee
but that would force us to work tangent to volume-preserving maps which would make the proof a little more complicated; we do not need to use that because the map
$$
\phi \mapsto \Lambda^{-1}{\mathbb H}_{imn}((\pa_p\phi^p)\tau_X^{mn})
$$
is bounded as a linear map from $C(0,T; [C^{1,\alpha}]^d)$ to itself.
Clearly, because $\tau_X$ is H\"{o}lder continuous, and because classical Calderon-Zygmund operators are bounded in H\"{o}lder spaces, there is no difficulty in bounding H\"{o}lder norms of derivatives of the expression $\Lambda^{-1}{\mathbb H}_{imn}((\pa_p\phi^p)\tau_X^{mn})$. Proving that the undifferentiated quantity is bounded is done using the fact that $\tau_X\in L^{1}\cap L^{\infty}$, and therefore $(\pa_p\phi^p)\tau_X\in L^{1}\cap L^{\infty}$. The operator
$\Lambda^{-1}{\mathbb H}$  maps continuously $L^1\cap L^{\infty}$ to $L^{\infty}$.
The operators ${\mathcal R}_{pimn}=\pa_p\Lambda^{-1}{\mathbb H}_{imn}$ are sums of classical Calderon-Zygmund operators and multiples of the identity. 
The commutators
\be
\phi\mapsto  \phi^p{\mathcal R}_{pimn}(\tau_X^{mn}) -{\mathcal R}_{pimn}(\phi_p\tau_X^{mn})
\la{cx}
\ee
are bounded as operators from  $C(0,T; [C^{1,\alpha}]^d)$ to itself. This is quite obvious for smooth $\tau_X$ but a little less obvious for $\tau_X\in \Sigma$.
Let us write the kernel of ${\mathcal R}_{pimn}$ as $K_{pimn}$, so the commutator is
\be
\int_{\Rr^d}K(x-y)(\phi(x)-\phi(y))\tau_X(y)dy
\la{commun}
\ee
where we did not write all the indices and the time dependence for ease of notation. The kernel $K$ is smooth away from the origin and is homogeneous of order $-d$.
Differentiating in some direction and writing $K'$ for the singular (of order $d+1$) kernel obtained by differentiating $K$, we have
\be
T[X]\phi = P.V.\int K'(x-y)(\phi(x)-\phi(y))\tau_X(y)dy
\la{TX}
\ee
plus a nice operator $(\nax \phi^p){\mathcal R}_{pimn}(\tau_X^{mn})$. This last operator is clearly bounded in $C^{\alpha}$, with bound controlled by $M$, so we concentrate our attention on $T[X]$. Now we write
\be
T[X]\phi = \int_0^1d\lambda\left [P.V.\int_{\Rr^d}(x-y)K'(x-y) \na \phi(x+\lambda(y-x))\tau_X(y)dy\right]
\la{txp}
\ee
The kernel $(x-y)K'(x-y)$ is homogeneous of order $-d$. It might have nonzero average on the unit sphere. Nevertheless, we subtract the value $\nax \phi(x)$:
\be
T_1[X]\phi = \int_0^1d\lambda\left\{P.V.\int_{\Rr^d}(x-y)K'(x-y)\left[\nax\phi(x+\lambda((y-x)))- \nax\phi(x)\right]\tau_X(y)dy\right\}
\la{tone}
\ee
The contributions left from the average on the unit sphere, if nonzero, are a constant multiple of $(\nax\phi (x))\tau_X(x)$ and $\nax\phi(x)T_2\tau_X(x)$,
\be
T_2(\tau_X)(x) = \int_{\Rr^d}(x-y)K'(x-y)(\tau_X(x)-\tau_X(y))dy
\la{ttwo}
\ee
both bounded with values in $C^{\alpha}$. The fact that $T_1[X]\phi$
is bounded in $C^{\alpha}$, and similarly, that $T_2(\tau_X)$ is a H\"{o}lder continuous function are classical. A proof can be found in (\cite{bc}). 
We have one more term in $DF[X]$, namely
\be
Y\mapsto k(\Lambda^{-1}{\mathbb H}((D\tau[X]Y)\circ X^{-1}))\circ X
\la{last}
\ee
Its boundedness is equivalent to the boundedness of the maps of the type
\be
\phi \mapsto \Lambda^{-1}{\mathbb H}(g_X \na\phi)
\la{lin}
\ee
in $C^{1,\alpha}$ where $g_X$ is in $\Sigma$. These are easily bounded because
when we take spatial derivatives we arrive at
\be
\phi \mapsto {\mathcal R}(g_X\na\phi)\la{linlast}
\ee
which are bounded in $C^{\alpha}$, and if we do not take derivatives, the $L^{\infty}$ boundedness follows as above from the fact that $g_X\na\phi \in L^1\cap L^{\infty}$.
We have therefore verified the fact that $DF[X]$ is bounded in ${\mathcal X}$
uniformly for $X\in {\mathcal D}$. The function $F$ is locally Lipschitz and because solutions of (\ref{lag}) start from the identity, they coincide for short time. The same argument does not need the initial data to be the identity, but rather the same invertible $C^{1,\alpha}$ transformation so uniqueness propagates because $F$ is locally Lipschitz. This concludes the proof of the theorem.

\section{A Regularization}
We consider here the system formed by the (\ref{lfp}) coupled with (\ref{stks}), (\ref{dxu}) via (\ref{sigmam}). In this section, we consider the same potential (\ref{um}) but we allow $R$ to be a function of $x$ and $t$ (but of course, not of $m$). The properties that $\sigma$ is symmetric, non-negative and the bound (\ref{indi}) remain valid. In order to see what dependences on $R$ would be allowed by the energy considerations, we repeat the calculation leading to (\ref{ens}): we take  (\ref{lfp}), multiply by $k(\log f + U)$,  integrate, take (\ref{stks}) multiply by $u$, integrate, and add the two: we deduce
\be
\ba
k\frac{d}{dt}\left[\int_{\Rr^d}\int_{\Rr^d}f\log fdm dx +\frac{1}{2}\int_{\Rr^d}Tr[\sigma(x,t)]dx \right] + k\int_{\Rr^d}\fr{2D_t R}{R}\left(\int_{\Rr^d} f U dm\right)dx\\
= -\int_{\Rr^d}|\nax u(x,t)|^2dx -k\e\int_{\Rr^d}\int_{\Rr^d}f\left| \nam(\log f+U)\right|^2dmdx 
\ea
\la{ensr}
\ee
It is thus clear that $D_tR\ge 0$ is energetically favorable. The solution on characteristics might be less explicit, however the equation obeyed by $\sigma$ is easily obtained by multiplying (\ref{lfp}) by $(m\otimes m)/R^2$ and integrating. The result is very similar to (\ref{sigeq}):
\be
D_t\sigma = (\nax u)\sigma + \sigma (\nax u)^T - \fr{2\e}{R^2}\sigma + \fr{2\e}{R^2}\rho {\mathbb I} - \fr{2D_t R}{R}\sigma
\la{sigreq}
\ee
We again see that $D_tR\ge 0$ has the effect of an additional damping. In fact, if $\fr{D_tR}{R}$ is a constant, then the effect is precisely one of enhanced damping, and that is similar to the situation covered previously in Theorem  {\ref{wsm}}, but in a better regime. More interesting, perhaps, is a damping that responds locally to very high rate of strain in the fluid.
Let consider a coupled system in which, in addition to (\ref{lfp} coupled to (\ref{stks}), (\ref{dxu}) via (\ref{um}), $R$ evolves according to
\be
D_t R = \delta\left(|\nax u(x,t)|\right )R
\la{req}
\ee
and $\delta (g)$ is a smooth nonnegative function of one nonnegative variable $g$, that vanishes for $g \le \frac{\kappa}{2}$
\be
\delta(g) = \left\{
\ba
0,\;\;\quad\quad\quad \quad \; \; \; \, {\mbox{if}}\; g\le \frac{\kappa}{2}\\
C_0\sqrt{\kappa^2 + g^2}, \quad {\mbox{if}}\; g\ge \kappa 
\ea
\right .
\la{deltag}
\ee
and satisfies
\be
|\delta'(g)| \le 2C_0
\la{deltaprime}
\ee
for all $g\ge 0$.  In particular  $\delta(|\nax u|)$ satisfies
\be
\delta = C_0 \sqrt{ \kappa^2 + |\nax u|^2}\la{delta}
\ee
if $|\nax u|\ge \kappa>0$. The constant $\kappa>0$ is fixed, and represents the order of magnitude of the largest permissible temporal growth rate. The constant $C_0>0$ is chosen in function of the dimension of space (this is needed because of the tensorial nature of the calculation) so that     
\be
3c|\nax u| -2\delta\le 0
\la{disi}
\ee
holds, if $|\nax u|\ge \kappa$, where $c$ is the norm equivalence constant that bounds $\max_{il}|\partial_l u^i| \le c |\nax u |$. Thus, $c$ is a norm equivalence constant in $\Rr^{d^2}$, and we can choose $C_0= \frac{3}{2}c$.

Because the problem evolves  in time, the initial data $f_0$ and $R_0$ need to be specified. We assume that $f_0$ is non-negative and smooth enough in $x$, decaying fast enough in $m$, ($f_0\in W^{1,p}(\Rr^d; L^1( (1+|m|^2)dm)$) and $R_0(a) \ge R_{min}>0$ is in $W^{1, p}(\Rr^d)$ with $p>d$. Then, from (\ref{sigreq}) and from the fact that $\sigma ^{ii}$ are nonnegative for each fixed $i$,
it follows that
\be 
D_t\,Tr(\sigma) \le 2c\kappa\, Tr(\sigma) + \fr{2d\e}{R^2}\rho
\la{trr}
\ee
which, after integration on characteristics,  results in
\be
 \sup_{x} Tr(\sigma(x,t)) \le e^{2c\kappa t}\left[\sup_{x} Tr(\sigma_0(x)) + \fr {d\e}{c\kappa R_{min}^2}\|\rho_0\|_{L^{\infty}(\Rr^d)}\right]= N_0(t).
\la{trs}
\ee
In view of (\ref{indi}), we thus have control of the $L^{\infty}$ norm of $\sigma$ in time. Once this is achieved, based on previous results for similar models, we may expect to prove regularity. The method of proof
of regularity given boundedness of $\sigma$ we employ here is the simplest and most explicit. The idea is to differentiate the equation, pay the price of differentiating the advective term, and do an Eulerian calculation in $L^p$ spaces, with $p>d$ . Working in Lebesgue spaces makes it easy to take advantage of incompressibility. The approach used in the proof of global existence for Smoluchowski equations on compact manifolds coupled with time-dependent Stokes equations in (\cite{c-smo}), and in the proof of  global existence    (\cite{cs1}) for Smoluchowski equations coupled to Navier-Stokes equations in $d=2$, can be adapted to the present situation as well. In other words, we can consider the time depending
$R$ evolving according to (\ref{req}), a Fokker-Planck equation (\ref{lfp}) with potential given in (\ref{um}) coupled via (\ref{sigmam}) to a fluid velocity evolving according to the time-dependent Stokes equation in $d=3$ or Navier-Stokes equation in $d=2$. The method of (\cite{chemm}), used in (\cite{cm}) to prove the global regularity for 
$d=2$ Smoluchowski equations on compact manifolds coupled with Navier-Stokes equations requires a little less smoothness on the Navier-Stokes initial data. That method is also Eulerian, but uses commutation and a penalty, that results in a controlled loss of regularity. The proof is a bit more technical as it employs Besov spaces and paradifferential calculus, but in principle it can be adapted to the noncompact particle phase space case with spatially depending $R$. We will not pursue these matters here, but rather content ourselves with the simplest proof, in the simplest nontrivial case.
Differentiating (\ref{sigreq}) we obtain
\be
\ba
D_t(\pa_k\sigma^{ij}) = -(\pa_ku^l)(\pa_l\sigma^{ij}) + (\pa_l u^i)(\pa_k\sigma^{lj}) + (\pa_l u^j)(\pa_k\sigma^{il}) -2\delta (\pa_k\sigma^{ij}) -\fr{2\e}{R^2}(\pa_k\sigma^{ij})\\ + (\pa^2_{kl}u^i)\sigma^{lj} + (\pa^2_{kl}u^j)\sigma^{il} -(\pa_k{\fr{2\e}{R^2}})\sigma^{ij} - 2\pa_k(\delta)(\sigma^{ij}) + \pa_k(\fr{2\e}{R^2}\rho \delta_{ij})
\ea
\la{difgr}
\ee
We mutiply by $\pa_k(\sigma^{ij})$ and sum. The terms involving explicitly first derivatives of $u$ are bounded using (\ref{disi}):
\be
\ba
\left[-(\pa_ku^l)(\pa_l\sigma^{ij}) + (\pa_l u^i)(\pa_k\sigma^{lj}) + (\pa_l u^j)(\pa_k\sigma^{il}) -2\delta (\pa_k\sigma^{ij}) -\fr{2\e}{R^2}(\pa_k\sigma^{ij})\right](\pa_k\sigma^{ij})\\
\le 3c\kappa |\nax\sigma|^2.
\ea
\la{first}
\ee
The terms involving explicitly second derivatives of $u$ are bounded using
(\ref{trs}):
\be
\left[(\pa^2_{kl}u^i)\sigma^{lj} + (\pa^2_{kl}u^j)\sigma^{il}\right](\pa_k\sigma^{ij}) \le CN_0(t)|\nax\nax u||\nax\sigma|.
\la{second}
\ee
The term containing a derivative of $\delta$ is bounded using (\ref{deltaprime}, which implies
\be
|\nax\delta(|\nax u|)| \le 2C_0|\nax\nax u|
\la{naxdelta}
\ee
and therefore 
\be
 -2\pa_k(\delta)(\sigma^{ij})(\pa_k\sigma^{ij})\le CN_0(t)|\nax\nax u||\nax\sigma|.
\la{third}
\ee 
Summarizing, we have so far, pointwise:
\be
\frac{1}{2}D_t|\nax \sigma|^2 \le 3c\kappa |\nax\sigma|^2 + CN_0(t)|\nax\nax u||\nax\sigma| + CN_0(t)\e \left|\nax\frac{1}{R^2}\right| |\nax\sigma| +
C\e\left |\nax\fr{\rho}{R^2}\right||\nax\sigma |.
\la{inco}
\ee
We have to bound the terms involving $\nax\rho$ and $\nax\fr{1}{R^2}$.
In view of the fact that $\rho (x,t) = \rho_0(A(x,t))$, we have
\be
\frac{\nax\rho}{R^2} = \left(\frac{\nax A}{R^2}\right)^T(\naa\rho_0)(A).
\la{naxrhoa}
\ee
Now, because of (\ref{naxaeq}) and (\ref{req}) we have
\be
D_t\left(\fr{\nax A}{R^2}\right ) = -\left(\fr{(\nax A)}{R^2}\right)\left((\nax u) +2\delta\mathbb I\right)
\la{naxar}
\ee
and, in view of (\ref{disi}) we deduce that
\be
\fr{|\nax A(x,t)|}{R(x,t)^2} \le \fr{e^{c\kappa t}}{R_{min}^2}
\la{naxarb}
\ee
and consequently
\be
\fr{|\nax\rho (x,t)|}{R(x,t)^2} \le C\fr{e^{c \kappa t}}{R_{min}^2}|\naa \rho_0(A(x,t))|\le C\fr{e^{c\kappa t}}{R_{min}^2}\|\naa\rho_0\|_{L^{\infty}(\Rr^d)}.
\la{naxrhob}
\ee
The term involving $\nax\left({\fr{1}{R^2}}\right)$ is treated using (\ref{req}) 
\be
D_t{\frac{\nax R}{R^3}} = \left (-(\nax u)^T - 2\delta\mathbb I\right)\left(\fr{\nax R}{R^3}\right) + \fr{\nax\delta}{R^2}
\la{naxrcube}
\ee  
and therefore, in view of (\ref{disi}) and (\ref{naxdelta}) we deduce
\be
D_t\left(\fr{|\nax R|}{R^3}\right) \le c\kappa \left(\fr{|\nax R|}{R^3}\right) +
C\frac{|\nax\nax u|}{R^2}
\la{naxrc}
\ee
and, integrating on characteristics, we obtain
\be
\fr{|\nax R(x,t)|}{R^3(x,t)}\le Ce^{c\kappa t}\left (\fr{|\naa R_0(A(x,t))|}{R_{min}^3} + \fr{1}{R_{min}^2}\int_0^t|\nax\nax u(X(A(x,t),s),s)|ds\right)
\la{naxrb}
\ee
Now we collect the terms in (\ref{inco}), divide by $|\nax\sigma|$ and use
(\ref{naxrhob}) and (\ref{naxrb}):
\be
\ba
D_t|\nax\sigma (x,t)| \le 3c\kappa|\nax\sigma(x,t)| + CN_0(t)|\nax\nax u(x,t)| +
C\fr{e^{c\kappa t}}{R_{min}^2}\e|\naa\rho_0(A(x,t))|\\
+C\e(N_0(t)+ \|\rho_0\|_{L^{\infty}(\Rr^d)})e^{c\kappa t}\left (\fr{|\naa R_0(A(x,t))|}{R_{min}^3} + \fr{1}{R_{min}^2}\int_0^t|\nax\nax u(X(A(x,t),s),s)|ds\right)
\ea
\la{near}
\ee
The inequality (\ref{near}) has the form
\be
D_t (y(x,t)) \le 3c\kappa y(x,t) + C(t)z(x,t) + D(t)\int_0^tz(X(A(x,t),s),s)ds +E(x,t)
\la{absnear}
\ee
where 
\be
\left\{
\ba
y(x,t) = |\nax\sigma (x,t)|,\\
z(x,t) = |\nax\nax u (x,t)| \\
C(t) = CN_0(t) = C e^{2c\kappa t}\left[\sup_{x} Tr(\sigma_0(x)) + \fr {d\e}{c\kappa R_{min}^2}\|\rho_0\|_{L^{\infty}(\Rr^d)}\right]\\
D(t) =C\fr{\e}{R_{min}^2}(N_0(t) +\|\rho_0\|_{L^{\infty}(\Rr^d)})e^{c\kappa t}, \\
E(x,t) = C\fr{\e}{R_{min}^2}e^{c\kappa t}\left \{|\naa\rho_0(A(x,t))| + \left (N_0(t)+\|\rho_0\|_{L^{\infty}(\Rr^d)}\right)\fr{|\naa R_0(A(x,t))|}{R_{min}}\right\}
\ea
\la{listsym}
\right.
\ee
The interested reader can check that (\ref{near}) is dimensionally balanced.
Now it is time to start looking at $L^p$ norms. We start by noting that
\be
\ba
\left \|\int_0^t z(X(A(x,t),s),s)ds\right\|_{L^p(dx)} \\
\le\int_0^t\|z(X(A(x,t),s),s)\|_{L^p(dx)}ds = \int_0^t\|z(X(a,s),s)\|_{L^p(da)}ds\\
=\int_0^t\|z(x,s)\|_{L^p(dx)}ds
\ea
\la{saus}
\ee
because of incompresssibility. 
We integrate (\ref{absnear}) on characteristics and take the $L^p$ norm. In order  simplify the answer we use the fact that $C(t), D(t)$ are non-decreasing functions of time. We obtain
\be
\|y(\cdot,t)\|_{L^p(\Rr^d)}\le e^{3c\kappa t}\left\{\|y_0\|_{L^p(\Rr^d)} + \|E(\cdot,t)\|_{L^p(\Rr^d)} + \int_0^t(C(s) + D(t)(t-s))\|z(\cdot, s)\|_{L^p(\Rr^d)}ds \right\}
\la{bou}
\ee
We recall (\ref{nas}), which implies
\be
\nax\nax u = k{\widetilde R}(\nax\ws).
\la{naxnaxu}
\ee
In view of the well-known boundedness of Riesz transforms in $L^p$ spaces, we deduce that
\be
\|z(\cdot,t)\|_{L^p(\Rr^d)} \le Ck\|y(\cdot,t)\|_{L^p(\Rr^d)}
\la{zyp}
\ee
We note also that $\|E(\cdot, t)\|_{L^p(\Rr^d)}$ is explicitly a sum of norms of initial data multiplied by exponentials of time.
Now a simple Gronwall argument provides an apriori bound for $y$ in $L^p$
\be
\|y(\cdot, t)\|_{L^p(\Rr^d)} \le F_p(t)
\la{yb}
\ee
with $F_p(t)$ an explicit function of time, with exponential growth, and depending only on norms of intitial data $\|\sigma_0\|_{L^{\infty}(\Rr^d)}$, 
$\|\rho_0\|_{L^{\infty}(\Rr^d)}$, $\|\nax\sigma_0\|_{L^p(\Rr^d)}$, and
$R_{min}$, $\|\nax R_0\|_{L^p(\Rr^d)}$.

\begin{thm} \la{elastic} Let $f$ solve (\ref{lfp}) with $U$ given by (\ref{um}) and
$R$ evolving according to (\ref{req}) with smooth $\delta$ satisfying (\ref{deltaprime}) and (\ref{disi}). Let $u$ be obtained by solving (\ref{stks}), (\ref{dxu}), with $\sigma$ given by (\ref{sigmam}). Assume that the initial distribution $f_0$ and $R_0$ satisfy
\be
\left \{
\ba
\sup_{x\in{\Rr^d}}\int_{\Rr^d} f_0(x,m)(1+|m|^2)dm <\infty\\
R_0(x)\ge R_{min}>0,\\
\int_{\Rr^d}|\nax R_0(x)|^pdx <\infty,\\
\int_{\Rr^d}\left (\int_{\Rr^d}(1+|m|^2)|\nax f_0(x,m)|dm\right)^pdx<\infty
\la{id}
\ea
\right.
\ee
with $p>d$. Then the solution exists for all time, is unique and obeys the a priori
bounds (\ref{trs}) and (\ref{yb}). 
\end{thm}
\noindent{\bf Proof.} The space $W^{1,p}(\Rr^d)$ is a space of local existence and uniqueness of solutions. The condition $\inf_{x}R(x,t)\ge R_{min}$ is invariant in time, because $R(x,t)\ge R_0(A(x,t))\ge R_{min}$. Then the apriori bounds are enough to finish the proof.

\begin{rmk} The a priori bounds hold in any $L^p$, $1<p<\infty$. 
\end{rmk}
\begin{rmk} Similar theorems hold if we replace steady-state Stokes equation with time-dependent Stokes equation in $d=2,3$ and with Navier-Stokes equation in $d=2$.  
\end{rmk}
\begin{rmk} Higher regularity of solutions can be obtained without difficulty from higher regularity of the initial data.
\end{rmk}
\section{Large data}
The problem of global regularity for arbitrary smooth initial data is open. The system formed by (\ref{lfp}) and (\ref{stks}) has potentially finite-time blowup. Indeed, consider the system in $d=1$
\be
\left\{
\ba
D_tf(x,m,t) + u_x\pa_m (mf(x,m,t)) =\e\pa_m(f\pa_m(\log f(x,m,t)+ U(m)))\\
U(m) = \fr{m^2}{2R^2},\\
\sigma (x,t) = \int_{-\infty}^{\infty} \fr{m^2}{R^2}f(x,m,t)dm\\
u_x = k H\sigma
\ea
\right.
\la{oned}
\ee
where $H$ is the Hilbert transform. 
This system is a 1-d analogue of (\ref{lfp}) with (\ref{um}), (\ref{sigmam}) and (\ref{nas}). This blows up in finite time. Indeed, we multiply the linear Fokker Planck equation by $\fr{m^2}{R^2}$ and integrate $dm$ to obtain the analogue of (\ref{sigeq}):
\be
D_t\sigma = 2u_x\sigma -\fr{2\e}{R^2}\sigma + \fr{2\e}{R^2}
\la{sigmaonedeq}
\ee
which then, in view of (\ref{oned}) is 
\be
D_t\sigma = 2k\sigma H\sigma -\fr{2\e}{R^2}\sigma +\fr{2\e}{R^2}
\la{baby}
\ee
and resembles the baby vorticity equation (\cite{clm}). The blow up
argument of (\cite{clm}) works here also, notwithstanding the fact that the present equation is computed on characteristics. We form
$z = H\sigma + i\sigma$ and use the fact that 
$$
H(\sigma H\sigma) = \fr{1}{2}\left((H\sigma)^2-\sigma^2\right)
$$
to deduce
$$
D_t z = kz^2 -\fr{2\e}{R^2}z + i\fr{2\e}{R^2}
$$
which blows up in finite time.

This simple example does not capture incompressibility and the tensorial nature of the problem, just like the baby vorticity equation of (\cite{clm}) does not capture them. The main quadratic nonlocal nonlinearity of (\ref{sigeq}) is modeled by $\sigma H\sigma$ and this is the only available model in one dimension.
Unlike the incompressible Euler equation however, the system we study has an additional dissipative structure, and we believe that this is of some significance and represents the main fallacy of the one-dimensional model. 
We describe below a simple scalar model that addresses this issue. We take henceforth $d=2$. We consider the variables
\be
\left\{
\ba
a(x,t) = \frac{1}{2}\left ( \sigma^{11}(x,t) - \sigma^{22}(x,t)\right),\\
b(x,t) = \sigma^{12}(x,t) = \sigma^{21}(x,t),\\
c(x,t) = \sigma^{11}(x,t) + \sigma^{22}(x,t) = Tr\,(\sigma(x,t))
\ea
\la{abc}
\right.
\ee
The gradient of velocity is represented by $\lambda(x,t)$, $\mu(x,t)$ and 
$\omega(x,t)$ given in (\ref{lam}), (\ref{mu}) and (\ref{ome}).
The equations (\ref{sigeq}) can be written as the system
\be
\left\{
\ba
D_t a = -\omega b + \lambda c - \fr{2\e}{R^2}a \\
D_t b  = \omega a + \mu c -\fr{2\e}{R^2} b,\\
D_t c = 4\lambda a + 4\mu b -{\fr{2\e}{R^2}}c + 4\fr{\e}{R^2}\rho\\
D_t\rho = 0
\ea
\right.
\la{abceq}
\ee
As we saw before, if we couple this with an equation for the velocity (steady or unsteady Stokes, or Navier-Stokes), then the regularity issue is decided by whether or not we can bound $c$ in $L^{\infty}(dx)$. We note in passing that the co-rotational system corresponds to $\lambda = \mu = 0$ in the system above, and the bound for $c$ follows immediately.
Let us multiply the $c$ equation by $\fr{c}{2}$, the $a$ equation by $2a$, the $b$ equation by $2b$ and subtract the last two from the first. We obtain
\be
D_t\left( \fr{c^2}{4} - a^2-b^2\right) = -\fr{4\e}{R^2}\left (\fr{c^2}{4}-a^2-b^2\right) + \fr{2\e}{R^2}\rho c
\la{deteq}
\ee
This cancellation of nonlinearity is not surprising because
\be
\fr{c^2}{4} - a^2-b^2 = Det\, (\sigma)
\la{detdet}
\ee
and the determinant is conserved along particle trajectories if $\e =0$.
The matrix $\sigma $ is symmetric and positive by construction, and is given in terms of $a,b,c$ by

\be
\sigma = \left (
\begin{array}{cc}
\fr{c}{2} + a & b\\
b & \fr{c}{2} -a
\end{array}
\right )
\la{sigmabc}
\ee
The positivity of the matrix is equivalent (in this case) to the positivity of the determinant, i.e. to
\be
\fr{c^2}{4}-a^2-b^2>0.
\la{detpos}
\ee
The two eigenvalues of $\sigma$
\be
z_{1,2} = \fr{c}{2} \pm \sqrt{a^2+b^2}
\la{ev}
\ee 
are both positive. We have of course $c =z_1+z_2>0$ and $z_1-z_2 = 2\sqrt{a^2+b^2}$. Because of (\ref{detpos}),
$c$ controls $\sqrt{a^2 +b^2}$. But, on the other hand, growth without bound of $c$ on any trajectory, cannot happen without growth without bound of $\sqrt{a^2+b^2}$ on the same trajectory. Indeed, if a particle path would be such that
$\sqrt{a^2 + b^2}$ is bounded on it, but $c$ grows without bound or blows up in finite time, then, for large enough time we would have
$$
-\fr{\e}{R^2}\left (\fr{c^2}{4}-a^2-b^2\right) + \fr{2\e}{R^2}\rho c\le 0
$$
on the particle path (because $\rho$ is bounded) and then, from (\ref{deteq}) we would arrive at the contradiction that $c$ remains bounded. 
From (\ref{abceq}) we can write 
\be
c(x,t) = e^{-\fr{2\e}{R^2}t}c_0(A(x,t)) + \int_0^te^{-\fr{2\e}{R^2}(t-s)}\left (\lambda a + \mu b -\fr{\e}{R^2}\rho_0\right)(X(A(x,t),s),s)ds
\la{cexpl}
\ee
Thus, exponential growth or blowup of $c$ is possible only if $(\lambda a + \mu b)$ has time integrals on particle paths that are positive and grow exponentially or stronger, without bound. 

In two dimensions we can express the velocity in terms of a stream function $\psi(x,t)$ and write $u^1=-\pa_2\psi$, $u^2=\pa_1\psi$. Then $\omega = \Delta \psi$, and therefore
\be
\lambda = \pa_1\pa_2(-\Delta)^{-1}\omega = B\omega\la{lamo}
\ee
and
\be
\mu = -\frac{1}{2}\left(\pa_1^2-\pa_2^2\right )(-\Delta)^{-1}\omega = -A\omega\la{muo}
\ee 
where
\be
A = \fr{1}{2}\left(\pa_1^2-\pa_2^2\right )(-\Delta)^{-1} =\fr{1}{2} (R_1^2-R_2^2)\la{A}
\ee
and  
\be
B = \pa_1\pa_2(-\Delta)^{-1} = R_1R_2
\la{B}
\ee
are bounded operators in $L^p(\Rr^2)$ spaces. They are also bounded, selfadjoint in $L^2(\Rr^2)$, they commute $AB=BA$, and each is given by a multiplier at Fourier level,
\be
\widehat{Ah}(\xi) = \fr{\xi_2^2-\xi_1^2}{2|\xi |^2}\widehat{h}(\xi)
\la{fa}
\ee
and
\be
\widehat{Bh}(\xi) = -\fr{\xi_1\xi_2}{|\xi|^2}\widehat{h}(\xi).
\la{fb}
\ee 
Note that
\be
4(A^2+B^2) = {\mathbb I}
\la{trig}
\ee
Let us consider now the time independent Stokes system (\ref{stks}). Taking the curl of (\ref{stks}) and inverting the Laplacian we obtain 
\be
\omega = 2k\left (Ab-Ba\right)
\la{omstks}
\ee
Consequently, from (\ref{lamo})
\be
\lambda = {2k}\left (-B^2 a + AB b\right)
\la{lab}
\ee
and
from (\ref{muo})
\be
\mu = {2k}\left(AB a - A^2b\right)
\la{muab}
\ee
It is convenient to measure time in units of $\fr{1}{2k}$ and then the system
(\ref{abceq}) coupled with (\ref{omstks}), (\ref{lamo}) and (\ref{muo}) is
\be
\left\{
\ba
D_t a = -(A(b)-B(a))b + [-B^2(a) +AB(b)]c -\varepsilon^{-1} a\\
D_t b = (A(b)-B(a))a  +[AB(a) - A^2(b)]c  -\varepsilon^{-1} b\\
D_t c = 4\left\{[-B^2(a)+AB(b)]a + [AB(a)-A^2(b)]b\right\} -\varepsilon^{-1} c + 2\varepsilon^{-1} \rho\\
D_t\rho = 0
\ea
\la{abcq}
\right.
\ee
where $\varepsilon$ is given by (\ref{varepsilon}). The all-important term $\lambda a + \mu b$ is given by
\be
\fr{1}{2k}(\lambda a + \mu b) = aB(\omega) -bA(\omega) =-aB(B(a)-A(b)) + bA(B(a)-A(b))
\la{gr}
\ee
This expression is quadratic, nonlocal, and has negative spatial integral.
While blow up requires the pointwise positivity of this expression (\ref{cexpl}), its spatial average is negative. Integrating the third equation in (\ref{abcq}) we obtain, using the selfadjointness of $A$ and $B$:
\be
\ba
\int_{\Rr^2}c(x,t)dx + 4\int_0^te^{-\fr{t-s}{\varepsilon}}\int_{\Rr^2}|B(a)(x,s)-A(b)(x,s)|^2dx \\
= e^{-\fr{t}{\varepsilon}}\int_{\Rr^2}c_0(x)dx +
2(1-e^{-\fr{t}{\varepsilon}})\int_{\Rr^2}\rho_0(x)dx
\ea
\la{trb}
\ee
This is just the energetic bound on $\int Tr\,\sigma dx$. We saw earlier that, if the size of initial data is of order one, then small enough $\varepsilon$ leads to global existence of solutions. The case of large initial data, moderate $\varepsilon$ is wide open. In order to clarify the issues, we will take $\varepsilon = \infty$. In this case the system is simplified somewhat, because, in view of
(\ref{deteq}) we have
\be
c = 2\sqrt{a^2+b^2 + d_0(x,t)}\la{cslave}
\ee
where 
\be
d_0(x,t) = Det\, (\sigma_0(A(x,t)))
\la{dxt}
\ee
and so, the first two equations of (\ref{abceq}) become just
\be
\left\{
\ba
D_t a = -b(Ab-Ba) + 2\sqrt{(a^2+b^2 +d_0)}\;B(Ab-Ba)\\
D_t b = a(Ab-Ba)   - 2\sqrt{(a^2+b^2 +d_0)}\;A(Ab-Ba).
\ea
\right .
\la{abeq}
\ee
and consequently
\be
D_t{\sqrt{(a^2+b^2 +d_0)}} = 2aB(Ab-Ba)-2bA(Ab-Ba)
\la{dtsqrt}
\ee
Integrating in space we again get the bound on the trace. Blowup for the system (\ref{abeq}) is still too difficult to analyse. The simplest possible didactic model of a scalar equation exhibiting this kind quadratic nonlocal structure with an  $L^1$ dissipation is a scalar equation of the form

\be
\pa_t \tau = -\tau (A^2\tau)
\la{mod}
\ee
where $A$ is a bounded selfadjoint operator given in Fourier representation by multiplication by a function that is homogeneous of degree zero, like above, 
\be
\widehat{Ah}(\xi) = \alpha(\xi)\widehat{h}(\xi)
\la{alpha}
\ee
with $\alpha(\lambda \xi) = \alpha(\xi)$, $\alpha(\xi)\in {\mathbb R}$,
$\int_{{\mathbb S}^2}\alpha(\xi)d\xi = 0$, $\alpha$ smooth on the unit sphere.
The unknown $\tau(x,t)$ (representing the trace of $\sigma$) is a positive scalar. We note the salient features of this model. If $\tau_0$ is smooth and positive, 
(let us consider for instance $\tau(\cdot,0)\in L^1\cap C^{0,s}$ with $s>0$) then the solution exists and is unique for some time, and remains positive as long as it exists.
If a bound on $\|\tau\|_{L^{\infty}}$ is given, then the higher regularity of the solution follows. Integrating the equation, we have the dissipative bound
\be
\int \tau(x,t)dx + \int_0^t\|A\tau\|_{L^2}^2dt \le \int\tau(x,0)dx
\la{diss}
\ee
and multiplying by $A^2\tau $ and integrating we have
\be 
\int |A\tau (x,t)|^2dx \le\int |A\tau(x,0)|^2dx
\la{l2}
\ee
We are going to give an example of global regularity for some large data for such an equation, in the spatially periodic case. For simplicity we take the period to be $2\pi$ in each direction. We write
\be
\tau(x,t) = \tau_0(t) + \sum_{k\in {\mathbb Z}^2\setminus\{0\}}\tau_k(t)e^{i k\cdot x}\la{taup}
\ee
The requirement that $\tau $ be real is implemented by $\overline{\tau_{k}}=\tau_{-k}$. We will consider the symbol $\alpha(k)$, and assume that it is real valued, even, $\alpha(-k)= \alpha(k)$, and bounded $|\alpha(k)| \le \Gamma$. We also assume  $\alpha(0) = 0$.
The equation (\ref{mod}) is the infinite system of ODEs
\be
\fr{d\tau_l}{dt} = - \sum_{k+j=l}\tau_k\alpha^2(j)\tau_j
\la{fode}
\ee
\begin{prop}\la{ode} Let $\tau(x,0) = \sum_k \tau_k(0)e^{ik\cdot x}$ satisfy
\be
\tau_{-k}(0) = \overline{\tau_k(0)},
\la{real}
\ee
\be
\sum_k(1+|k|)^s|\tau_k(0)|\le C_s(0)<\infty
\la{smid}
\ee
for some $s>0$ and
\be
\tau_0(0) \ge \sum_{k\neq 0}|\tau_k(0)|.
\la{posid}
\ee
Then the solution of (\ref{fode}) exists for all time, and obeys
\be
\tau_{-k}(t) = \overline{\tau_{k}(t)}
\la{realt}
\ee
\be
\sum_k(1+|k|)^s|\tau_k(t)|\le C_s(0)e^{2^{s+1}\tau_0(0)\Gamma^2 t}<\infty
\la{smt}
\ee
and
\be
\tau_0(t) \ge\sum_{k\neq 0}|\tau_k(t)|
\la{inv}
\ee
which of course implies that $\tau (x,t)$ remains smooth and positive.
\end{prop}
\noindent{\bf{Proof.}} We start with the dissipation, the equation at $l=0$, which reads 
\be
\fr{d\tau_0}{dt} =- \sum_{k\neq 0}\alpha(k)^2|\tau_k(t)|^2
\la{dsp}
\ee
Now for $l\neq 0$ we have
$$
\fr{d\tau_l}{dt} = -\sum_{k+j=l,\; k\neq 0,\; j\neq 0}\tau_k\alpha^2(j)\tau_j - \tau_0\alpha^2(l)\tau_l
$$
and therefore 
\be
\fr{d|\tau_l|}{dt} \le -\tau_0\alpha^2(l)|\tau_l| + \sum_{k+j=l, \; k\neq 0, \; j\neq 0}|\tau_k|\alpha^2(j)|\tau_j|
\la{taul}
\ee
Summing in $l\neq 0$ we obtain
$$
\ba
\fr{d}{dt}\sum_{l\neq 0} |\tau_l| \le -\tau_0\sum_{l\neq 0}\alpha^2(l)|\tau_l| + \sum_{l\neq 0}\sum_{k+j=l,\; k\neq 0, \; j\neq 0}|\tau_k|\alpha^2(j)|\tau_j|\\
=  -\tau_0\sum_{l\neq 0}\alpha^2(l)|\tau_l| + \sum_{j\neq 0}\sum_{k\neq 0}\alpha^2(j)|\tau_j||\tau_k| - \sum_{j\neq 0}\alpha^2(j)|\tau_{-j}||\tau_j|\\
=\left(\sum_{j\neq 0}\alpha^2(j)|\tau_j|\right)\left (-\tau_0 +\sum_{k\neq 0}|\tau_k|\right) + \fr{d}{dt}\tau_0
\ea
$$
This results in the inequality
\be
\fr{d}{dt}\left (\sum_{l\neq 0}|\tau_l| -\tau_0\right)\le \left(\sum_{j\neq 0}\alpha^2(j)|\tau_j|\right)\left(\sum_{l\neq 0}|\tau_l| - \tau_0\right)
\la{ineqt}
\ee
If the initial data is non-positive, the quantity remains non positive. Thus we have the invariance of this cone in function space. Once this is achieved, we know
\be
\sum_{l\neq 0}|\tau_l(t)| \le \tau_0(t)\le \tau_0(0)
\la{linf}
\ee
which implies an $L^{\infty}$ bound on $\tau$. We know this should be sufficient for regularity. In the present situation, the proof of regularity is quite straightforward: we multiply (\ref{taul}) by $(1 + |l|)^s$ and sum. Using $|l|\le 2 \max(|k|, |j|)$ 
and $|\alpha(j)|\le \Gamma$ we obtain
$$
\fr{d}{dt}\sum_{l\neq 0} (1 + |l|)^s|\tau_l|\le 2^s\Gamma^2\left(\sum_j|\tau_j|\right)\left(\sum_{l\neq 0}(1+|l|)^s|\tau_l|\right)
$$
Using (\ref{linf}) we deduce
\be
\fr{d}{dt}\sum_{l\neq 0} (1 + |l|)^s|\tau_l| \le 2^{s+1}\tau_0(0)\Gamma^2\left(\sum_{l\neq 0}(1+|l|)^s|\tau_l|\right)
\la{ns}
\ee
which finishes the proof.

{\bf Acknowledgments:} PC's research was partially supported by NSF-DMS grant  0804380.

\end{document}